\documentclass[10pt,twoside]{article}

\usepackage[centertags]{amsmath}
\usepackage{amsfonts}
\usepackage{amssymb}
\usepackage{amsthm}
\usepackage{epsfig}
\usepackage{color}
\allowdisplaybreaks[4]

\definecolor{c20}{rgb}{0.,0.7,0.}
\definecolor{c30}{rgb}{0.,0.,1.}
\definecolor{c40}{rgb}{1,0.1,0.7}
\definecolor{c50}{rgb}{1,0,0}
\definecolor{c60}{rgb}{1,0.9,0.1}

\def\iE#1{\textcolor{c20}{#1}}
\def\iE#1{#1}

\def\he#1{\textcolor{c40}{#1}}
\def\he#1{#1}

\def\Quan#1{\textcolor{c30}{#1}}
\def\Quan#1{#1}
\def\aa#1{\textcolor{c40}{#1}}
\def\aa#1{#1}

\def\cy#1{\textcolor{c40}{#1}}
\def\cy#1{#1}

\def\AA#1{\textcolor{c30}{#1}}
\def\cE#1{\textcolor{c30}{#1}}
\def\cH#1{\textcolor{c30}{#1}}
\def\ccE#1{\textcolor{c30}{#1}}
\def\cT#1{\textcolor{c50}{#1}}
\def\ccT#1{\textcolor{c50}{#1}}
\def\ccH#1{\textcolor{c30}{#1}}
\def\cZ#1{\textcolor{c50}{#1}}
\def\ccZ#1{\textcolor{c50}{#1}}
\def\Tzq#1{\textcolor{c50}{#1}}
\def\Zhong#1{\textcolor{c50}{#1}}
\def\TT#1{\textcolor{c50}{#1}}
\def\cE#1{#1}
\def\cT#1{#1}
\def\ccT#1{#1}
\def\ccE#1{#1}
\def\cH#1{#1}
\def\ccH#1{#1}
\def\cZ#1{#1}

\def\AA#1{#1}
\def\ccZ#1{#1}

\def\Tzq#1{#1}
\def\Zhong#1{#1}
\def\TT#1{#1}
\def\Zhong#1{#1}

\newcommand{\COM}[1]{}

\newcommand{\ABs}[1]{ \biggl \lvert #1 \biggr \rvert}
\newcommand{\R}{\mathbb{R}}
\newcommand{\N}{\mathbb{N}}
\newcommand{\inr}{\in \R}
\newcommand{\BQN}{\begin{eqnarray}}
\newcommand{\EQN}{\end{eqnarray}}
\newcommand{\BQNY}{\begin{eqnarray*}}
\newcommand{\EQNY}{\end{eqnarray*}}
\def\IF{\infty}
\newcommand{\kb}[1]{\boldsymbol{#1}}
\newcommand{\vk}[1]{\cH{\kb{#1}}}
\def\MTD{M^\delta(T)}
\def\MT{M(T)}

\def\Phidz{\AA{d\Phi_p(\vk{z})}}

\begin{document}

\thispagestyle{empty}

\centerline{\large \bf \cH{On \Tzq{Piterbarg's} Max-discretisation
Theorem for \Tzq{Multivariate Stationary Gaussian Processes}}}

\mbox{}
       \vskip 0.4 cm

\centerline{Zhongquan Tan\footnote{College of Mathematics Physics and Information Engineering, Jiaxing University, Jiaxing 314001, PR China} and Enkelejd Hashorva\footnote{corresponding author: Department of Actuarial Science, Faculty of Business and Economics, University of Lausanne, Switzerland}}

       \vskip 0.4 cm
\centerline{Jiaxing University \&  University of Lausanne, Switzerland}
\centerline{\today{}}

       \vskip 0.4 cm

{\bf Abstract:}  \cE{Let $\{X(t), t\geq0\}$ be a stationary Gaussian
process with zero-mean and unit variance. A deep result derived in
Piterbarg (2004), which we refer to as \Tzq{Piterbarg's}
max-discretisation theorem gives the joint asymptotic behaviour
($T\to \infty$) of the continuous time  maximum $M(T)=\max_{t\in [0,T]} X(t),
$ and the maximum $M^{\delta}(T)=\max_{t\in \mathfrak{R}(\delta)}X(t), $ with $\mathfrak{R}(\delta)
\subset [0,T]$ a uniform grid of points of distance
$\delta=\delta(T)$. Under some asymptotic restrictions on the
correlation function \Tzq{Piterbarg's} max-discretisation theorem
shows that for the limit result it is important to know the speed
$\delta(T)$ approaches 0 as $T\to \IF$. The present contribution
derives the aforementioned theorem for \Tzq{multivariate stationary
Gaussian processes}\COM{ and presents an application which is of
certain importance for simulations}.}



\section{Introduction}

Let $\{X(t), t\geq0\}$ be a \cE{standard (zero-mean, unit-variance)} stationary Gaussian process with correlation
function $r(\cdot)$ and continuous sample \ccE{paths}.  \cE{A tractable and very large class of correlation functions satisfy}
\begin{eqnarray}
\label{eq1.1}
r(t)=1-C|t|^{\alpha}+o(|t|^{\alpha})\ \ \mbox{as}\ \ t\rightarrow 0
\end{eqnarray}
for some positive constant $C$ and $\alpha\in(0, 2]$, \ccH{see e.g.,
Piterbarg (1996)}. \cE{If further, the Berman condition (see Berman
(1964) \ccH{or} Berman (1992))}
\begin{eqnarray}
\label{eq1.2}
\lim_{\ccH{T}\rightarrow\infty} r(T) \ln   T=0
\end{eqnarray}
holds, then it is well-known, see e.g., Leadbetter et al. (1983)
that the maximum  $M(T)= \max_{ t\in [0,T]} X(t)$ obeys the Gumbel
law as $T\to \IF$, namely \BQN\label{Gumbel}
\cE{\lim_{T\rightarrow\infty}\sup_{x\inr}}\ABs{P\{a_{T}(M(T)-b_{T})\leq
x\}- \Lambda(x)}=0 \EQN \cE{is valid with
$\Lambda(x)=\exp(-\exp(-x)),x\inr$ \Quan{the cumulative distribution
function of} \iE{a Gumbel random variable} and normalising constants
defined for all large $T$ by} \BQN\label{aT} a_{T}=\sqrt{2 \ln   T},
\ \ b_{T}=a_{T}+a_T^{-1} \ln
\Bigl((2\pi)^{-1/2}C^{1/\alpha}H_{\alpha}a_T^{-1+2/\alpha}\Bigr).
\EQN Here $H_{\alpha}$ denotes \cE{the well-known} Pickands constant
\cE{given by the limit relation}
$$H_{\alpha}=\lim_{S\rightarrow\infty}\cE{S^{-1}}\cH{\mathbb{E}\biggl\{\exp\left(\max_{t\in[0,S]}
\iE{\Bigl(}\sqrt{2}B_{\alpha/2}(t)-t^{\alpha}
\iE{\Bigr)}\right)\biggr\}} \in (0,\infty),$$ with
 $B_{\ccH{\alpha}}$ a standard fractional Brownian motion \Tzq{with Hurst index $\alpha$}, \cy{see e.g.,
 Mishura and Valkeila (2011)} for recent characterisations of $B_\alpha$. For the main properties of Pickands and related constants, see for example Adler (1990), Piterbarg (1996), D\c{e}bicki (2002),  D\c{e}bicki et al.\ (2003), Wu (2007), D\c{e}bicki and Kisowski (2009), \AA{D\c{e}bicki and  Tabi\'{s} (2011) and Hashorva et al.\ (2013a).} \cE{We note in passing that the first correct
proof of Pickands theorem where $H_\alpha$ appears (see Pickands
(1969)) is derived in Piterbarg (1972).}

\Tzq{We say that $X$ is weakly dependent if its
correlation function satisfies the Berman condition
(\ref{eq1.2}). A natural generalisation of (2) is
the following assumption
\begin{eqnarray}
\label{eqStrong}
\lim_{\ccH{T}\rightarrow\infty} r(T) \ln   T=r\in(0,\infty)
\end{eqnarray}
in which case we say that $X$ is a strongly dependent Gaussian
process. Mittal and Ylvisaker (1975) prove the limit theorem for the \aa{normalised} maximum of strongly dependent stationary Gaussian
processes showing that the limiting distribution is a mixture of
Gumbel and Gaussian distribution. In fact, a similar result is shown
therein also for the extreme case that \eqref{eqStrong} holds with
$r=\IF$ with the limiting distribution being Gaussian.} \Zhong{For
other related results on extremes of strongly dependent Gaussian
sequences and processes, we refer to McCormick and Qi (2000), James
et al. (2007), Tan and  Wang (2012), Hashorva and Weng (2013), Hashorva et al.\ (2013b) and the references therein.}

\Quan{In this paper the random variable  $M(T)=\sup_{0\leq t\leq T} X(t),\iE{T>0}$ \aa{denotes} the
continuous-time maximum and  $M^{\delta}(T) =
\max_{t\in\delta\aa{\mathbb{N}}\cap[0,T ]}X(t)$ \aa{stands for} the maximum over
the uniform grid $\delta\aa{\mathbb{N}}\cap[0,T ]$.}  \COM{we refer to
$\MT$ as the continuous time maximum; the discrete time maximum
denoted by $\MTD$ will be defined using a uniform grid of points
$\mathfrak{R}(\delta)=\{ k\delta: k=0,1,2,\cdots\} \subset [0,T]$.
Specifically, for fixed $T>0$ and some positive $\delta=\delta(T)$
we define the maximum of the Gaussian process $\{X(t),t\ge 0\}$ over
this grid by $\MTD\ccH{=}\max_ {t\in \mathfrak{R}(\delta)} X(t)$.}
Under the assumption (1) we need to distinguish between three types
of grids: A uniform grid of points $\mathfrak{R}(\delta)=\aa{\delta\aa{\mathbb{N}}}$ is called
sparse if $\delta=\delta(T)$ is such that \BQN \label{grid}
\cE{\lim_{T\to \infty}\delta(T)}(2 \ln   T)^{1/\alpha}
=D, \EQN with $D=\infty$. When \eqref{grid} holds for some $D\in(0,
\infty)$, then the  grid is referred to as Pickands grid, whereas
when \eqref{grid} holds with $D = 0$, it is called a dense
\Quan{grid}. \COM{Note in passing that only the constant $\alpha \in
(0,2]$ is used in the definition of the different grids.} \Quan{Throughout \iE{this paper} we assume \iE{that} $\alpha\in (0, 2]$.}

\Quan{ Piterbarg (2004) \ccE{derived}} the joint asymptotic
behaviour of $\MTD$ and
 $\MT$ for weakly dependent stationary Gaussian processes. As shown therein, after a suitable normalisation
 (as in  \eqref{Gumbel})  $\MTD$ and $\MT$ \cE{are} asymptotically independent, dependent or totally dependent
  if the grid is a sparse, a Pickands or a dense grid, respectively.
 \cE{We shall refer to \ccE{that} result  as \Tzq{Piterbarg's} max-discretisation theorem}.

For a \cE{large} class of locally stationary Gaussian processes H\"{u}sler (2004) proved a similar result \ccE{to Piterbarg (2004)} considering \ccE{only} sparse and dense grids.
In another investigation concerning the storage process with fractional Brownian motion as input, it was shown in
 H\"{u}sler and Piterbarg (2004) that the continuous time maximum and the discrete time maximum over \Quan{the dense grid} are asymptotically completely dependent.
 \Zhong{Tan and Tang (2012) and} Tan and Wang (2013)  recently \Quan{proved} \Tzq{Piterbarg's} max-discretisation theorem for strongly dependent
stationary Gaussian processes, whereas Tan and Hashorva (2012)
derives \Quan{similar results for sparse and dense grids} for
standardised maximum of stationary Gaussian processes. \iE{Novel and deep
results concerning stationary \Quan{non-Gaussian} processes are derived in Turkman (2012).}

\cE{As noted in Piterbarg (2004) derivation of the joint asymptotic behaviour of
$\MTD$ and $\MT$ is important for theoretical problems and at the
same time is crucial for applications, see
Piterbarg (2004), H\"{u}sler (2004) and Tan and Hashorva (2012) for more details. }

\Zhong{The main contribution of this paper is the derivation of
\Tzq{Piterbarg's} max-discretisation theorem for \Tzq{multivariate
stationary Gaussian processes}.} Our results show that, despite the
high technical difficulties, it is possible to state Piterbarg's
result in multidimensional settings allowing for asymptotic
conditions \Quan{and the two maxima are no longer asymptotically
independent}.

\cE{Brief organisation of the paper:} \ccE{ In Section 2 we present
the principal results. \Tzq{Section 3 presents some auxiliary
results followed by  Section 4 which is dedicated} to  the proofs of
\iE{the our main theorems}. \iE{Several technical lemmas \aa{and} the proof of Lemma 3.1 are displayed in Appendix.}}

\section{Main Results}
\ccE{Consider} $(X_{1}(t),\cdots,X_{p}(t)), \cE{p\in \N}$  a $p$-dimensional centered
Gaussian \cE{vector} process \cE{with} covariance functions
$r_{kk}(\tau)=Cov(X_{k}(t),X_{k}(t+\tau)),\cE{k\le p}$.
Hereafter we shall assume that \cE{the components} have continuous sample paths and further  $Cov(X_{k}(t),\cZ{X_{l}}(t+\tau))$ does not dependent on $t$ so we shall write
$$r_{kl}(\tau)=Cov(X_{k}(t),X_{l}(t+\tau))$$
for the cross-covariance function. Further we shall suppose that each component $X_i$ has a unit variance function; in short we shall refer to such vector processes as standard stationary Gaussian vector process.
\cE{Similarly to \eqref{eq1.1} we suppose that} for all indices $k\le p$
\begin{eqnarray}
\label{eq1.3}
r_{kk}(t)=1-\cT{C}|t|^{\alpha}+o(|t|^{\alpha})\ \ \mbox{as}\ \ t\rightarrow 0,
\end{eqnarray}
\cE{with \ccH{some} positive constants \cT{$C$}}, and further
\begin{eqnarray}
\label{eq1.4}
\lim_{T\to \IF} r_{kl}(T) \ln   T= r_{kl}\in\Quan{(0,\infty)}
\end{eqnarray}
holds for $1\leq k,l\leq p$. \ccE{In order to} exclude the possibility that $X_{k}(t)=\pm X_{l}(t+t_{0})$ for some $k\neq l$, and some choice of $t_{0}$ and $+$ or $-$, we assume that
\begin{eqnarray}
\label{eq1.5}
\max_{k\neq l}\sup_{\AA{t \in [0,\IF)}}|r_{kl}(t)|<1.
\end{eqnarray}
For any $k\le p$ and a given \aa{uniform} grid of points $\mathfrak{R}(\delta)$ we define the componentwise maximum (in continuous and discrete time) by
$$(M_{k}(T), M_{k}^{\delta}(T)):=\Bigl(\max_{t\in[0,T]}X_{k}(t), \max_{t\in\mathfrak{R}(\delta)\cE{\cap}[0,T]}X_{k}(t)\Bigr).$$

\ccH{Let $\ccH{\vk{Z}}=(Z_1, \ldots, Z_p)$  be a  $p$-dimensional centered Gaussian random vector with covariances}
$$Cov(Z_{k},Z_{l})=\frac{r_{kl}}{\sqrt{r_{kk}r_{ll}}}.$$
 \iE{Further, let $\Psi$ denote the survival function of a
$N(0,1)$ random variable} and put $\vk{x}:=(x_1, \ldots, x_p),
\vk{y}:=(y_1, \ldots, y_p).$

In our theorem below we consider the case of sparse grids, followed then  by two results on Pickands and dense grids.

\def\xyrp{\cH{\vk{x}\in \R^p, \vk{y} \in \R^p}}
\newcommand{\EE}[1]{\mathbb{E}\Bigl\{#1 \Bigr\}}

\textbf{Theorem 2.1}. {\sl \ccE{Let $(X_{1}(t),\cdots,X_{p}(t))$ be standard stationary Gaussian vector process with}  covariance functions satisfying (\ref{eq1.3}), (\ref{eq1.4}) and (\ref{eq1.5}). \COM{\ccH{If further $\vk{Z}$ has a positive-definite covariance matrix},}
\AA{If further $\vk{Z}$ has a positive-definite covariance matrix,} then for any sparse grid $\ccH{\mathfrak{R}(\delta)}$
\begin{eqnarray}
\label{eq2.1}
\lim_{T\to \IF}\sup_{ \xyrp } \ABs{
P\left\{a_{T}(M_{k}(T)-b_{T})\leq x_{k}, a_{T}(M_{k}^{\delta}(T)- b_{T}^{\delta})\leq y_{k}, k=1,\cdots,p\right\}
-\EE{\exp\left(-f(\vk{x}, \vk{y}, \vk{Z})\right) }
}=0,
\end{eqnarray}
\ccE{where \ccH{$a_T$ is defined in \eqref{aT}},
$$f(\vk{x}, \vk{y},\vk{Z})=\ccH{\sum_{k=1}^{p}}
\left(e^{-x_{k}-r_{kk}+\sqrt{2r_{kk}}Z_{k}}+e^{-y_{k}-r_{kk}+\sqrt{2r_{kk}}Z_{k}}\right)$$
and} \BQN\label{btd} b_{T}^{\delta}=a_T+a_T^{-1}   \ln
\bigl((2\pi)^{-1/2}\delta^{-1}a_T^{-1}\bigr). \EQN}


\textbf{Corollary 2.1}. {\sl \ccE{Under the assumptions of Theorem 2.1} we have
\begin{eqnarray}
\label{eq2.4}
\lim_{T\to \IF}\sup_{ \ccH{\vk{x}\in \R^p}} \ABs{
P\left\{a_{T}(M_{k}(T)-b_{T})\leq x_{k}, k=1,\cdots,p\right\}
-\EE{\exp\left(-h(\vk{x}, \vk{x}, \vk{Z})\right) }}=0,
\end{eqnarray}
where
\BQN\label{HhH}
h(\vk{x}, \vk{y}, \vk{Z})=\sum_{k=1}^{p}e^{-\min(x_{k},\ccH{y_k})-r_{kk}+\sqrt{2r_{kk}}Z_{k}}
\EQN
and $a_T, b_T$ are defined in \eqref{aT}.
}

\ccE{Before presenting the result for Pickands grids}, \ccE{we
introduce} the following constants which can be found in Leadbetter
et al. \ccH{(1983)}. For any $d>0,\lambda>0$, \Quan{$k\in
\mathbb{Z}$} and $x,y\inr$ define
$$H_{d,\alpha}(\lambda)=\mathbb{E}\biggl \{ \exp\left(\max_{kd\in[0,\lambda]}\iE{\Bigl(}\sqrt{2}B_{\alpha/2}(kd)-(kd)^{\alpha}\iE{\Bigr)}\right) \biggr\}$$
and
$$H_{\ccZ{d},\alpha}^{\he{x,y}}(\lambda)= \lim_{\lambda\to \IF}\frac{1}{\lambda}\int_{\he{s\inr}}e^{s}P\left(
\max_{t\in[0,\lambda]}\Bigl(\sqrt{2}B_{\alpha/2}(t)-t^{\alpha}\Bigr)>s+\he{x},
\max_{k:kd\in [0,\lambda]}\Bigl(\sqrt{2}B_{\alpha/2}(k\ccZ{d})-(k\ccZ{d})^{\alpha}\Bigr)>s+\he{y} \right)ds.$$

\ccE{In view of} Leadbetter et al. (1983) both constants $H_{d,\alpha}$  and $H_{d,\alpha}^{x,y}$ defined as
$$H_{d,\alpha}=\lim_{\lambda\rightarrow\infty}\frac{H_{d,\alpha}(\lambda)}{\lambda} \quad
\text{  and } \quad
H_{\ccZ{d},\alpha}^{x,y}= \lim_{\lambda\to \IF}\frac{H_{d,\alpha}^{x,y}(\lambda)}{\lambda}$$
are finite and positive.

\textbf{Theorem 2.2}. {\sl \ccE{Let $(X_{1}(t),\cdots,X_{p}(t))$ satisfy the assumptions of Theorem 2.1, and let
 $a_T$ be as in \eqref{aT}}. For any Pickands grid $\mathfrak{R}(d a_T^{-2/\alpha})$ with $d>0$ we have
\BQN
\label{eq2.2}
\lim_{T\to \IF}\sup_{ \xyrp } \ABs{
P\left\{a_{T}(M_{k}(T)-b_{T})\leq x_{k}, a_{T}(M_{k}^{\delta}(T)- \cZ{b_{d,T}})\leq y_{k}, k=1,\cdots,p\right\}
-\EE{\exp\left(-g(\vk{x}, \vk{y}, \vk{Z})\right) }}=0,
\EQN
\ccE{where}
$$g(\vk{x}, \vk{y}, \vk{Z})=\ccH{\sum_{k=1}^{p}}\left(e^{-x_{k}-r_{kk}+\sqrt{2r_{kk}}Z_{k}}+e^{-y_{k}-r_{kk}+\sqrt{2r_{kk}}Z_{k}}-H_{\ccZ{d},\alpha}
^{ \he{\ln   H_{\alpha}+x_{k}, \ln   H_{\ccZ{d},\alpha}+y_{k}}}e^{-r_{kk}+\sqrt{2r_{kk}}Z_{k}}\right),$$
with
\BQN\label{btd2}
b_{d,T}=a_T+ a_T  \ln  \Bigl((2\pi)^{-1/2}C^{1/\alpha}H_{d,\alpha}a_T^{-1+2/\alpha}\Bigr).
\EQN
}


\textbf{Theorem 2.3}. {\sl \ccE{Under the assumptions of Theorem 2.1} for any dense grid $\mathfrak{R}(\delta)$
\begin{eqnarray}
\label{eq2.3}
\lim_{T\to \IF}\sup_{ \xyrp } \ABs{
P\left\{a_{T}(M_{k}(T)-b_{T})\leq x_{k}, a_{T}(M_{k}^{\delta}(T)- b_{T})\leq y_{k}, k=1,\cdots,p\right\}
-\EE{\exp\left(-h(\vk{x}, \vk{y}, \vk{Z})\right) }}=0,
\end{eqnarray}
\ccE{with $a_T, b_T$ as defined in \eqref{aT}} and \aa{$h$ defined in \eqref{HhH}.}
}

\def\ldot{,\ldots,}
\cT{
\textbf{Remark 2.1}.  \AA{a)} In condition (\ref{eq1.3}) we can use different $C's$ and $\alpha's$, i.e., condition (\ref{eq1.3}) can be replaced by
\begin{eqnarray*}
r_{kk}(t)=1-C_{k}|t|^{\alpha_{k}}+o(|t|^{\alpha_{k}})\ \ \mbox{as}\ \ t\rightarrow 0.
\end{eqnarray*}
In that case, the above results still hold with some obvious \ccE{modifications} of $b_{T}$, $b_{T}^{\delta},b_{d,T}$
and the grid $\mathfrak{R}(\delta)$.\\
\AA{b) If $\vk{Z}$ has a singular covariance matrix, then we still
can derive the above results. To see that consider the simple case
$Z_p=c_1Z_1+  \cdots+ c_{p-1}Z_{p-1}$ with $c_{i},i\le p-1$ some
constants, and $\vk{Z}^*=(Z_1 \ldot Z_{p-1})$ has a non-singular
covariance matrix. In the proofs below we need to condition on
$Z_1=z_1 \ldot Z_{p-1}=z_{p-1}$, and then put $c_1z_1+  \cdots+
c_{p-1}z_{p-1}$ instead of $z_p$ therein. } }

\COM{\textbf{Remark 2.2}. a) \TT{It is worth  pointing out that the two
maxima are  no longer asymptotically independent since $r_{kk}>0$
for all $k\le p$ in Theorem 2.1. In fact, Theorem 2.1 still holds
for the case that $r_{kk}=0$ for all $k$ whereas the two maxima are
independent. In that case we note that $r_{kl}$'s play no role in the
result. For the proof of this case, we only need to suppose that
$\vk{Z}=\mathbf{0}$.
\\
b) Theorems 2.2 and 2.3 also hold for the case that $r_{kk}=0$ for
all $k$.}
}
\def\SL{{\cal S}_l}
\def\SI{{\cal S}_i}
\def\SJ{{\cal S}_j}

\def\RL{{\cal R}_l}
\def\RJ{{\cal R}_j}
\newcommand{\pk}[1]{P\left\{#1\right\} }
\newcommand{\abs}[1]{\lvert #1 \rvert}
\newcommand{\Es}[1]{\mathbb{E}\{#1 \}}

\section{\cE{Auxiliary \iE{R}esults}}
\cE{In this section we present several lemmas needed for the proof
of the main results, \Tzq{where Lemma 3.1 plays a crucial role.} In
order to establish \Tzq{Piterbarg's} max-discretisation theorem for
 \Tzq{standard stationary vector Gaussian processes} we need to closely
follow the steps of the proofs in Piterbarg (2004), and of course to
strongly rely on the deep ideas and techniques presented in
Piterbarg (1996). First, for $1\leq k,l\leq p$ define
$$\rho_{kl}(T)=r_{kl}/\ln T.$$ \COM{and let
$$\varpi_{kl}(t)=\max\{|r_{kl}(t)|, \rho_{kl}(T), |r_{kl}(t)+(1-r_{kl}(t))\rho_{kl}(T)|\},$$ $$\vartheta_{kl}(t)=\sup_{t<s\leq T}\{\varpi_{kl}(s)\}\ \ \mbox{and}\ \ \theta_{kl}(t)=\sup_{t<s\leq T}\{|r_{kl}(s)|\}.$$
For $1\leq k,l\leq p$ \ccH{we have}
$$\theta_{kl}(t)\cZ{\leq}\vartheta_{kl}(t)<1$$
for each $t\in [0,T)$ and all sufficiently large $T$. Further, for $1\leq k,l\leq p$, let \cZ{$a,b$ be constants such} that
$$0<b<a<\frac{ 1-\vartheta_{kl}(0)}{1+\vartheta_{kl}(0)}\cZ{\leq} \frac{1-\theta_{kl}(0)}{1+\theta_{kl}(0)}<1$$
 for all sufficiently large $T$.}
Following the former reference, we divide the} interval $[0,T]$ onto
intervals \cE{of} length $S$ alternating with shorter intervals
\ccE{of} length $R$. Let $a> b$ be constants which \iE{will be
determined in} \Quan{the proof of Lemma 3.1.} We shall denote
throughout in the sequel
$$S=T^{a}, \quad  R=T^{b}, \quad \iE{T>0}.$$
 Denote the long intervals by $ \SL $, $l=1,2,\cdots,n_{T}=[T/(S+R)]$, and the short intervals by $ \RL $,
$l=1,2,\cdots,n_{T}$.  It will be seen \cE{from the proofs,} that a
possible remaining interval with length different than $S$ or $R$
plays no role in our \cE{asymptotic considerations}; we call also
this interval a short interval. \cE{Define further} $\mathbf{S}=\cup_{\iE{l=1}}^{n_T}
\SL , \mathbf{R}=\cup_{\iE{l=1}}^{n_T}  \RL $ \iE{and thus}  $[0,T]=\mathbf{S}\cup
\mathbf{R}$.

\Tzq{Our proofs also rely on the main ideas of Mittal and Ylvisaker
(1975) by constructing new Gaussian processes to approximate the
original Gaussian processes.} For each index $k \le p$ we define a
new Gaussian process $\eta_k$ by taking  $\{Y_{k}^{(j)}(t), t\geq
0\}$, $j=1,2,\cdots, n_{T}$ independent copies of $\{X_{k}(t),
t\geq0\}$ and setting  $\eta_{k}(t)=Y_{k}^{(j)}(t)$ for $t\in
\RJ\cup\SJ=[(j-1)(S+R),j(S+R))$. We construct the processes so that
$\eta_k,k=1, \cdots, p$ are independent by taking $Y_k^{\iE{(j)}}$ to
be independent for any $j$ and $k$ two possible indices. The
independence of two different processes  $\eta_k$ and $\eta_l$
implies
$$\gamma_{kl}(s,t):=\EE{\eta_{k}(s)\eta_{l}(t)}=0,\ \ k\neq l,$$
whereas for any fixed $k$
\[
  \gamma_{kk}(s,t):=\EE{\eta_{k}(s)\eta_{k}(t)}=\left\{
 \begin{array}{cc}
  { \EE{Y_k^{\iE{(i)}}(t),Y_k^{\iE{(i)}}(s)}=  r_{kk}(\Quan{s,t})},   & \text{ if } t,s\in \mathcal{R}_{i}\cup \mathcal{S}_{i}, \text{ for some } i \le n_{T};\\
  \EE{Y_k^{\iE{(i)}}(t),Y_k^{\iE{(j)}}(s)}={0},    & \text{ if } t\in \mathcal{R}_{i}\cup \mathcal{S}_{i}, s\in \mathcal{R}_{j}\cup \mathcal{S}_{j},
  \text{ for some } i\neq j \le n_{T}.
 \end{array}
  \right.
\]\\
For $k=1,2,\cdots,p$ define
$$\xi_{k}^{T}(t)=\big(1-\rho_{kk}(T)\big)^{1/2}\eta_{k}(t)+\rho^{1/2}_{kk}(T)Z_{k}, \ \ 0\leq t\leq T,$$
where \ccH{$\vk{Z}=(Z_1, \ldots, Z_p)$} is a $p$-dimensional
centered Gaussian random vector \Quan{defined in Section
2},\COM{with covariances
$$Cov(Z_{k},Z_{l})=\frac{r_{kl}}{\sqrt{r_{kk}r_{ll}}},$$}
which is independent of $\{\eta_{k}(t), t\geq 0\}$,
$k=1,2,\cdots,p$. \COM{It can be checked that  $\{\xi_{k}^{T}(t),
0\leq t\leq T, k=1,2,\cdots,p\}$ are stationary Gaussian vector
processes.} Denote by $\{\varrho_{kl}(s,t),1\leq k,l\leq p\}$ the
covariance functions of $\{\xi_{k}^{T}(t), 0\leq t\leq
T,k=1,2,\cdots,p\}$. We have
$$\varrho_{kl}(s,t)=\EE{\xi_{k}^{T}(s)\xi_{l}^{T}(t)}=\rho_{kl}(T),\ \ k\neq l$$
and
\[
  \varrho_{kk}(s,t)=\left\{
 \begin{array}{cc}
  {r_{kk}(\Quan{s,t})+(1-r_{kk}(\Quan{s,t}))\rho_{kk}(T)} ,    &t\in \mathcal{R}_{i}\cup \mathcal{S}_{i}, s\in \mathcal{R}_{j}\cup \mathcal{S}_{j}, i= j;\\
  {\rho_{kk}(T)},    & t\in \mathcal{R}_{i}\cup \mathcal{S}_{i}, s\in \mathcal{R}_{j}\cup \mathcal{S}_{j}, i\neq j.
 \end{array}
  \right.
\]
\COM{Since both $\gamma_{kk}(s,t)$ and $\varrho_{kk}(s,t)$ are
functions of $|s-t|$  we shall write below simply
$\gamma_{kk}(|s-t|)$ and $\varrho_{kk}(|s-t|)$, respectively.}

For \cE{notational simplicity} we write for any $x_k,y_k\in \R$
$$u_{k}^{(1)}= b_{T}+x_{k}/a_{T}, \quad u_{k}^{(2)}= b_{T}^{'}+y_{k}/a_{T},$$
where $b_{T}^{'}=b_{T}^{\delta}$ for \ccH{a} sparse grid and $b_{T}^{'}=b_{d,T}$ for Pickands grid. \cE{Further,} for any $\varepsilon>0$ set
$$q=\frac{\varepsilon}{( \ln   T)^{1/\alpha}}.$$

\Tzq{We give first a crucial result which shows that the maximums of
the original Gaussian processes $(X_{1}(t),\cdots,X_{p}(t))$ can be
approximated by that of the Gaussian processes $\{\xi_{k}^{T}(t),
0\leq t\leq T, k=1,2,\cdots,p\}$. \iE{The proof of \Quan{the} next
Lemma, due to its complications, is relegated to the Appendix.}}

\textbf{Lemma 3.1}. {\sl \Quan{Suppose that the grid
$\mathfrak{R}(\delta)$ is a sparse grid or a Pickands grid.} For any
$B>0$ for all $x_{k},y_{k}\in[-B,B], \ccH{k\le p}$ \ccH{we have}
\begin{eqnarray*}
&& \lim_{T \to \IF} \bigg|P\left\{\max_{t\in \mathfrak{R}(q)\cap\mathbf{S}}X_{k}(t)\leq u_{k}^{(1)}, \max_{t\in\mathfrak{R}(\delta)\cap\mathbf{S}}X_{k}(t)\leq u_{k}^{(2)}, k=1,\cdots,p\right\}\\
&&\ \ \ \ \ \ \ \ \ \ \  -P\left\{\max_{t\in \mathfrak{R}(q)\cap \mathbf{S}}\xi_{k}^{T}(t)\leq u_{k}^{(1)}, \max_{t\in\mathfrak{R}(\delta)\cap \mathbf{S}}\xi_{k}^{T}(t)\leq u_{k}^{(2)}, k=1,\cdots,p\right\}\bigg|
=0.
\end{eqnarray*}
}

In order to deal with our multivariate framework in Lemmas 3.2 and
3.3 below we present the multivariate versions of Lemmas 6 and 4 in
Piterbarg (2004), respectively. Lemma 3.4 is a new result.

\textbf{Lemma 3.2}. {\sl \Quan{Suppose that the grid
$\mathfrak{R}(\delta)$ is a sparse grid or a Pickands grid.} For any
$B>0$ \ccE{there exits a positive} constant \cZ{$K$} such that for
all $x_{k},y_{k}\in[-B,B]$ \ccE{we have}
\begin{eqnarray*}
&&\bigg|P\left\{M_{k}(T)\leq u_{k}^{(1)}, M_{k}^{\delta}(T)\leq u_{k}^{(2)}, k=1,\cdots,p\right\} \\
&&\ \ \ \ \ -P\left\{\max_{t\in \mathbf{S}}X_{k}(t)\leq u_{k}^{(1)}, \max_{t\in\mathfrak{R}(\delta)\cap\mathbf{S}}X_{k}(t) \leq u_{k}^{(2)}, k=1,\cdots,p\right\}\bigg| \leq \cZ{K}( \ln   T)^{1/\alpha-1/2}T^{b-a},
\end{eqnarray*}
\ccE{with $0< b< a < 1$ given constants \ccE{and all $T$ large.}} }

\textbf{Proof:} \ccE{In order to obtain the upper bound, we shall
use the following inequality}
\begin{eqnarray}
\label{eq3.1}
&&\bigg|P\left\{M_{k}(T)\leq u_{k}^{(1)}, M_{k}^{\delta}(T)\leq u_{k}^{(2)}, k=1,\cdots,p\right\} \nonumber\\
&&\ \ \ \ \ -P\left\{\max_{t\in \mathbf{S}}X_{k}(t)\leq u_{k}^{(1)}, \max_{t\in\mathfrak{R}(\delta)\cap\mathbf{S}}X_{k}(t) \leq u_{k}^{(2)}, k=1,\cdots,p\right\}\bigg|\nonumber\\
&&\ \ \ \ \ \leq \sum_{k=1}^{p}P\left\{\max_{t\in \mathbf{R}}X_{k}(t)>b_{T}+\cZ{x_{k}}/a_{T}\right\}
+\sum_{k=1}^{p}P\left\{\max_{t\in \mathbf{R}\cap \mathfrak{R}(\delta)}X_{k}(t)>b_{T}^{'}+\cZ{y_{k}}/a_{T}\right\}
\end{eqnarray}
\ccE{valid for any $x_i\inr,y_i\inr,i\le k.$} By Pickands theorem
$$P\left\{\max_{t\in \mathbf{R}}X_{k}(t)>b_{T}+x_{k}/a_{T}\right\}=O(mes(\mathbf{R})(b_{T}+\cZ{x_{k}}/a_{T})^{2/\alpha})\Psi(b_{T}+\cZ{x_{k}}/a_{T})$$
as $T\rightarrow\infty$, where $mes(\mathbf{R})$ denotes the Lebesgue measure of $\mathbf{R}$.
\ccE{In the light of}  (11) and (16) of Piterbarg
(2004) for a sparse grid and Pickands grid, respectively, we get the
same order for the second probability in the \ccE{right-hand side}
of (\ref{eq3.1}), hence the proof is complete. \cE{\hfill $\Box$}

\textbf{Lemma 3.3}. {\sl \Quan{Suppose that the grid
$\mathfrak{R}(\delta)$ is a sparse grid or a Pickands grid.} For any
$B>0$ for all $x_{k},y_{k}\in[-B,B], k\le p$ and for the Pickands
grid $\mathfrak{R}(q)=\mathfrak{R}(\varepsilon/(\ln T)^{1/\alpha})$ we have
\begin{eqnarray*}
&&\bigg|P\left\{\max_{t\in \mathbf{S}}X_{k}(t)\leq u_{k}^{(1)}, \max_{t\in\mathfrak{R}(\delta)\cap\mathbf{S}}X_{k}(t)\leq u_{k}^{(2)}, k=1,\cdots,p\right\}\\
&&\ \ \ \ \ -P\left\{\max_{t\in \mathfrak{R}(q)\cap\mathbf{S}}X_{k}(t)\leq u_{k}^{(1)}, \max_{t\in\mathfrak{R}(\delta)\cap\mathbf{S}}X_{k}(t)\leq u_{k}^{(2)}, k=1,\cdots,p\right\}\bigg|
\rightarrow 0
\end{eqnarray*}
as $\varepsilon\downarrow 0$. }

\textbf{Proof:} \ccH{In view of} Lemma 4 of Piterbarg (2004)
\begin{eqnarray*}
&&\bigg|P\left\{\max_{t\in \mathbf{S}}X_{k}(t)\leq u_{k}^{(1)}, \max_{t\in\mathfrak{R}(\delta)\cap\mathbf{S}}X_{k}(t)\leq u_{k}^{(2)}, k=1,\cdots,p\right\}\\
&&\ \ \ \ \ \ \ \ \ \ \ \ \ \ \ \ \ -P\left\{\max_{t\in \mathfrak{R}(q)\cap\mathbf{S}}X_{k}(t)\leq u_{k}^{(1)}, \max_{t\in\mathfrak{R}(\delta)\cap\mathbf{S}}X_{k}(t)\leq u_{k}^{(2)}, k=1,\cdots,p\right\}\bigg|\\
&&\leq \sum_{k=1}^{p}
\left[P\left\{\max_{t\in \mathfrak{R}(q)\cap\mathbf{S}}X_{k}(t)\leq u_{k}^{(1)}\right\}-P\left\{\max_{t\in \mathbf{S}}X_{k}(t)\leq u_{k}^{(1)}\right\}\right]\\
&&\leq \Quan{g(\varepsilon)} \sum_{k=1}^{p}n_{T}T^{a}(u_{k}^{(1)})^{2/\alpha}\Psi(u_{k}^{(1)})\leq \cZ{K}\Quan{g(\varepsilon)},
\end{eqnarray*}
\Quan{where $g(\varepsilon)\rightarrow 0$ as $\varepsilon\downarrow
0$,} hence the claim follows. \cE{\hfill $\Box$}

\textbf{Lemma 3.4}. {\sl \Quan{Suppose that the grid
$\mathfrak{R}(\delta)$ is a sparse grid or a Pickands grid.} \cZ{For
any $B>0$ for all $x_{k},y_{k}\in[-B,B], \ccH{k\le p}$  and for the
Pickands grid $\mathfrak{R}(q)=\mathfrak{R}(\varepsilon/(\ln
T)^{1/\alpha})$} we have
\begin{eqnarray*}
&&\bigg|
P\left\{\max_{t\in \mathfrak{R}(q)\cap \mathbf{S}}\xi_{k}^{T}(t)\leq u_{k}^{(1)}, \max_{t\in\mathfrak{R}(\delta)\cap \mathbf{S}}\xi_{k}^{T}(t)\leq u_{k}^{(2)}, k=1,\cdots,p\right\}
\\
&&\ \ \ \ \ \ - \int_{\mathbf{z} \in \mathbb{R}^p}\prod_{i=1}^{n_{T}}P\left\{\max_{t\in \SI}\eta_{k}(t)\leq u_{k}^{*}, \max_{t\in\mathfrak{R}(\delta)\cap \SI}\eta_{k}(t)\leq u_{k}^{*'}, k=1,\cdots,p\right\}
\Phidz 
\bigg|\rightarrow 0
\end{eqnarray*}
as $\varepsilon\downarrow 0$, where \AA{ $\Phi_p$ is the distribution function of $\vk{Z}$}
\begin{eqnarray}
\label{eq3.4.2}
u_{k}^{*}:=\frac{b_{T}+x_{k}/a_{T}-\rho^{1/2}_{kk}(T)z_{k}}{(1-\rho_{kk}(T))^{1/2}}=\frac{x_{k}+r_{kk}-\sqrt{2r_{kk}}z_{k}}{a_{T}}+b_{T}+o(a_{T}^{-1}),
\end{eqnarray}
\begin{eqnarray}
\label{eq3.4.3}
u_{k}^{*'}:=\frac{b_{T}^{'}+y_{k}/a_{T}-\rho^{1/2}_{kk}(T)z_{k}}{(1-\rho_{kk}(T))^{1/2}}=\frac{y_{k}+r_{kk}-\sqrt{2r_{kk}}z_{k}}{a_{T}}+b_{T}^{'}+o(a_{T}^{-1})
\end{eqnarray}
 and  $b_{T}^{'}=b_{T}^{\delta}$ for a sparse grid and
$b_{T}^{'}=b_{a,T}$ for a Pickands grid. }

\textbf{Proof:} First, by the definition of $\xi_{k}^{T}$ and
$\eta_k$ we have
\begin{eqnarray}
\label{eq3.4.1}
&&P\left\{\max_{t\in \mathfrak{R}(q)\cap \mathbf{S}}\xi_{k}^{T}(t)\leq u_{k}^{(1)}, \max_{t\in\mathfrak{R}(\delta)\cap \mathbf{S}}\xi_{k}^{T}(t)\leq u_{k}^{(2)}, k=1,\cdots,p\right\}\nonumber\\
&&\ \ = \int_{\mathbf{z} \in \mathbb{R}^p}P\left\{\max_{t\in \mathfrak{R}(q)\cap\mathbf{S}}\eta_{k}(t)\leq u_{k}^{*}, \max_{t\in\mathfrak{R}(\delta)\cap\mathbf{S}}\eta_{k}(t)\leq u_{k}^{*'}, k=1,\cdots,p \right\}\Phidz\notag\\
&&\ \ = \int_{\mathbf{z} \in \mathbb{R}^p}\prod_{i=1}^{n_{T}}P\left\{\max_{t\in \mathfrak{R}(q)\cap \SI}\eta_{k}(t)\leq u_{k}^{*}, \max_{t\in\mathfrak{R}(\delta)\cap \SI}\eta_{k}(t)\leq u_{k}^{*'}, k=1,\cdots,p \right\}\Phidz.
\end{eqnarray}
As for the discrete case, see \AA{p.}\ 137 on Leadbetter et al.
(1983) direct calculations lead to (18) and (19).
\ccE{Next, similarly to the proof of} Lemma 3.3, for all large
$T$
\begin{eqnarray*}
&&\bigg|\prod_{i=1}^{n_{T}}P\left\{\max_{t\in \mathfrak{R}(q)\cap \SI}\eta_{k}(t)\leq u_{k}^{*}, \max_{t\in\mathfrak{R}(\delta)\cap \SI}\eta_{k}(t)\leq u_{k}^{*'}, k=1,\cdots,p\right\}\\
&&\ \ \ \ \ \ \ \ -\prod_{i=1}^{n_{T}}P\left\{\max_{t\in \SI}\eta_{k}(t)\leq u_{k}^{*}, \max_{t\in\mathfrak{R}(\delta)\cap \SI}\eta_{k}(t)\leq u_{k}^{*'}, k=1,\cdots,p\right\}\bigg|\\
&&\leq \sum_{i=1}^{n_{T}}\left[P\left\{\max_{t\in \mathfrak{R}(q)\cap \SI}\eta_{k}(t)\leq u_{k}^{*}, k=1,\cdots,p\right\}-P\left\{\max_{t\in \SI}\eta_{k}(t)\leq u_{k}^{*'}, k=1,\cdots,p\right\}\right]\\
&&\leq \Quan{g(\varepsilon)} T^{a}\sum_{i=1}^{n_{T}}(u_{k}^{*})^{2/\alpha}\Psi(u_{k}^{*})
\leq K\Quan{g(\varepsilon)}
\end{eqnarray*}
holds for some constant $K$, thus the claim follows by applying the dominated convergence theorem and letting $\varepsilon\downarrow0$.
\cE{\hfill $\Box$}

\def\AI{{\cal A}_i}
\def\AJ{{\cal A}_j}
\def\AK{{\cal A}_k}
\def\AL{{\cal A}_l}
\def\APK{{\cal A}_{p+k}}
\def\APL{{\cal A}_{p+l}}

\section{Proofs}

\textbf{Proof of Theorem 2.1}. From Lemmas 3.1-3.4 and the dominated
convergence theorem, we known that in order to prove Theorem 2.1, it
\ccE{suffice to show} that
\begin{eqnarray}
\label{eq2.1.1}
&&\lim_{T\rightarrow\infty}\bigg|\prod_{i=1}^{n_{T}}P\left\{\max_{t\in  \SI }\eta_{k}(t)\leq u_{k}^{*}, \max_{t\in\mathfrak{R}(\delta)\cap  \SI }\eta_{k}(t)\leq u_{k}^{*'}, k=1,\cdots,p\right\} \nonumber\\
&&\ \ \ \ \ \ \ \ \ \ \ \ \ \ \ \ \ \ \ \ \ -\exp\left(-\sum_{k=1}^{p}(e^{-x_{k}-r_{kk}+\sqrt{2r_{kk}}z_{k}}+e^{-y_{k}-r_{kk}+\sqrt{2r_{kk}}z_{k}})\right)\bigg|=0 .
\end{eqnarray}
Define next the events \cE{$${\cal A}_i=\Bigl\{ \max_{t\in
[0,S]}\eta_{i}(t)> u_{i}^{*} \Bigr\}, \quad  {\cal A}_{p+i}=
\Bigl\{\max_{t\in\mathfrak{R}(\delta)\cap [0,S]}\eta_{i}(t)>
u_{i}^{*'}\Bigr\}, \quad  i=1,\cdots,p.$$} Using the stationarity of
$\{\eta_{k}(t), k=1,2,\cdots,p\}$ (we write $\AI^c$ for the
complimentary event of $\AI$) \BQNY
\prod_{i=1}^{n_{T}}P\left\{\max_{t\in  \SI }\eta_{k}(t)\leq
u_{k}^{*}, \max_{t\in\mathfrak{R}(\delta)\cap  \SI }\eta_{k}(t)\leq
u_{k}^{*'}, k=1,\cdots,p\right\}
 &=& (P\{ \cap_{i=1}^{2p} \AI^c \})^{n_{T}}\\
&=&\exp\big(n_{T} \ln  (  P\{ \cap_{i=1}^{2p} \AI^c \})\big)\\
&=&\exp\big(-n_{T} P\{ \cup_{i=1}^{2p} \AI \}+W_{n_{T}}\big).
\end{eqnarray*}
Since $\lim_{T \to \IF} P\{ \cap_{i=1}^{2p} \AI \} =1$ we get that
the remainder $W_{n_{T}}$ satisfies
$$W_{n_{T}}=o(n_{T} P\{ \cup_{i=1}^{2p} \AI \}), \quad T \to \IF.$$
Next, by  Bonferroni inequality
\BQN\label{bonfi}
\sum_{i=1}^{2p} P\{ \AI \}  &\ge & P\{ \cup_{i=1}^{2p} \AI \} 
 \ge  \sum_{i=1}^{2p} P\{ \AI \} - \sum_{1 \le i < j \le 2p} P\{ \AI, \AJ \} \notag\\
& = & \sum_{i=1}^{2p} P\{ \AI \} - \sum_{1 \le k < l \le p} P\{ \AK, \AL \} - \sum_{1 \le k < l \le p} P\{ \APK, \APL \} -
2 \sum_{1 \le k < l \le p} P\{ \AK, \APL \}\notag\\
&&- \sum_{i=1}^p  P\{ \AI,{\cal A}_{p+i}
\}=:A_{1}-A_{2}-A_{3}-2A_{4}-A_{5}. \EQN

\cE{Further, Lemma 2 of  Piterbarg (2004) and (\ref{eq3.4.2}),
(\ref{eq3.4.3}) imply}
\begin{eqnarray*}
\cE{A_{1}}&\thicksim & \sum_{k=1}^{p}ST^{-1}(e^{-x_{k}-r_{kk}+\sqrt{2r_{kk}}z_{k}}+e^{-y_{k}-r_{kk}++\sqrt{2r_{kk}}z_{k}})\\
&=&\sum_{k=1}^{p}\frac{T^{a}}{T}(e^{-x_{k}-r_{kk}+\sqrt{2r_{kk}}z_{k}}+e^{-y_{k}-r_{kk}+\sqrt{2r_{kk}}z_{k}}), \quad T\rightarrow\infty.
\end{eqnarray*}
For $A_{2}$, by the independence of $\eta_{k}(t)$ and $\eta_{l}(t)$,
$k\neq l$, Lemma 2 of  Piterbarg (2004) and (\ref{eq3.4.2}),
(\ref{eq3.4.3}),
 we have
\begin{eqnarray*}
A_{2}&=&\cZ{\sum_{1\leq k<l\leq p}P\left\{\max_{t\in [0,S]}\eta_{k}(t)> u_{k}^{*}, \max_{t\in [0,S]}\eta_{l}(t)> u_{l}^{*}\right\}}\nonumber\\
&=&\sum_{1\leq k<l\leq p}P\left\{\max_{t\in [0,S]}\eta_{k}(t)> u_{k}^{*}\right\}P\left\{ \max_{t\in [0,S]}\eta_{l}(t)> u_{l}^{*}\right\}\nonumber\\
&\thicksim&\sum_{1\leq k<l\leq p}ST^{-1}e^{-x_{k}-r_{kk}+\sqrt{2r_{kk}}z_{k}}ST^{-1}e^{-x_{l}-r_{ll}+\sqrt{2r_{ll}}z_{l}}
=o(A_{1}).
\end{eqnarray*}
\COM{\cE{Utilising Lemma 4 of  Piterbarg (2004) and recalling that
$\mathfrak{R}(q)$ is a Pickands grid, for any $\epsilon>0$ we have}
\begin{eqnarray}
\label{eq2.1.2}
&&\bigg|P\left\{\max_{t\in [0,S]}\eta_{k}(t)>u_{k}^{*}, \max_{t\in [0,S]}\eta_{l}(t)>u_{l}^{*}\right\}
-P\left\{\max_{t\in \mathfrak{R}(q)\cap[0,S]}\eta_{k}(t)>u_{k}^{*}, \max_{t\in \mathfrak{R}(q)\cap[0,S]}\eta_{l}(t)>u_{l}^{*}\right\}\bigg|\nonumber\\
&&\leq \bigg|P\left\{\max_{t\in [0,S]}\eta_{k}(t)>u_{k}^{*}\right\}-P\left\{\max_{t\in \mathfrak{R}(q)\cap[0,S]}\eta_{k}(t)>u_{k}^{*}\right\}\bigg|\nonumber\\
&&\ \ \ \ \ \ \ +\bigg|P\left\{\max_{t\in [0,S]}\eta_{l}(t)>u_{l}^{*}\right\}-P\left\{\max_{t\in \mathfrak{R}(q)\cap[0,S]}\eta_{l}(t)>u_{l}^{*}\right\}\bigg|\nonumber\\
&&\leq\epsilon T^{a}(u_{k}^{*})^{2/\alpha}\Psi(u_{k}^{*})+\epsilon T^{a}(u_{l}^{*})^{2/\alpha}\Psi(u_{l}^{*})\leq \epsilon\frac{T^{a}}{T}.
\end{eqnarray}
\cE{Consequently,} from (\ref{eq2.1.2}) and Lemma 2 of Piterbarg (2004) \cE{we obtain}
\begin{eqnarray*}
A_{3}&=&\cZ{\sum_{1\leq k<l\leq p}P\left\{\max_{t\in [0,S]}\eta_{k}(t)> u_{k}^{*}, \max_{t\in [0,S]}\eta_{l}(t)> u_{l}^{*}\right\}(1+o(\frac{T^{a}}{T}))}\nonumber\\
&=&\sum_{1\leq k<l\leq p}P\left\{\max_{t\in [0,S]}\eta_{k}(t)> u_{k}^{*}\right\}P\left\{ \max_{t\in [0,S]}\eta_{l}(t)> u_{l}^{*}\right\}(1+o(\frac{T^{a}}{T}))\nonumber\\
&=&\sum_{1\leq k<l\leq p}ST^{-1}e^{-x_{k}-r_{kk}+\sqrt{2r_{kk}}z_{k}}ST^{-1}e^{-x_{l}-r_{ll}+\sqrt{2r_{ll}}z_{l}}(1+o(\frac{T^{a}}{T}))
=o(A_{1}).
\end{eqnarray*}}
Note that $\mathfrak{R}(\delta)$ is a sparse grid, similar arguments
as for $A_{2}$ imply
$$A_{k}=o(A_{1}), \ \ k=2,3.$$
\COM{Now, Lemma 5 in Piterbarg (2004)  yields
\begin{eqnarray*}
\lim_{T\rightarrow\infty}\bigg|P\left\{\max_{t\in[0,S]}\eta_{k}(t)>u_{k}^{*}, \max_{t\in \mathfrak{R}(\delta)\cap[0,S]}\eta_{l}(t)>u_{l}^{*'}\right\} -P\left\{\max_{t\in \mathfrak{R}(q)\cap[0,S]}\eta_{k}(t)>u_{k}^{*}, \max_{t\in \mathfrak{R}(\delta)\cap[0,S]}\eta_{l}(t)>u_{l}^{*'}\right\}\bigg|= 0.
\end{eqnarray*}
Thus, a similar argument as for $A_{2}$ leads to $A_{4}=o(A_{1}).$}
\cE{Further,} Lemma 2 of  Piterbarg (2004) implies $A_{5}=o(A_{1}).$
Consequently, as $T\to \IF$
$$n_{T} P\{ \cup_{i=1}^{2p} \AI \} \thicksim \sum_{k=1}^{p}\left(e^{-x_{k}-r_{kk}+\sqrt{2r_{kk}}z_{k}}+e^{-y_{k}-r_{kk}+\sqrt{2r_{kk}}z_{k}}\right),$$
which completes the proof of (\ref{eq2.1.1}). \cE{\hfill $\Box$}

\textbf{Proof of Theorem 2.2}. \ccE{In view of Lemmas 3.1-3.4 and the dominated convergence theorem in order to establish the proof we need to show}
\begin{eqnarray}
\label{eq2.2.1}
&&\bigg|\prod_{i=1}^{n_{T}}P\left\{\max_{t\in \SI}\eta_{k}(t)\leq u_{k}^{*}, \max_{t\in\mathfrak{R}(\delta)\cap \SI}\eta_{k}(t)\leq u_{k}^{*'}, k=1,\cdots,p\right\}\nonumber\\
&&\ \ \ \ \ \ \ \ -\exp\left(-\sum_{k=1}^{p}\left(e^{-x_{k}-r_{kk}+\sqrt{2r_{kk}}z_{k}}+e^{-y_{k}-r_{kk}+\sqrt{2r_{kk}}z_{k}}-H_{d,\alpha}^{\Quan{\ln   H_{\alpha}+x_{k},\ln   H_{d,\alpha}+y_{k}}}e^{-r_{kk}+\sqrt{2r_{kk}}z_{k}}\right)\right)\bigg|\rightarrow 0
\end{eqnarray}
as $T\rightarrow\infty$.
\ccE{We proceed as for the proof of (\ref{eq2.1.1}) using the lower bound \eqref{bonfi}; \ccH{we have thus}}
\ccT{\BQN\label{bonfin}
P\{ \cup_{i=1}^{2p} \AI \}
& = & \sum_{i=1}^{2p} P\{ \AI \} - \sum_{1 \le k < l \le p} P\{ \AK, \AL \} - \sum_{1 \le k < l \le p} P\{ \APK, \APL \} -
2 \sum_{1 \le k < l \le p} P\{ \AK, \APL \}\notag\\
&&- \sum_{i=1}^p  P\{ \AI,{\cal A}_{p+i} \}+\sum_{1 \le k < \cdots<
l \le p}=:A_{1}-A_{2}-A_{3}-2A_{4}-A_{5}+A_{6}. \EQN} \COM{We have
\begin{eqnarray*}
&&1-P\left\{\max_{t\in [0,S]}X_{k}(t)\leq u_{k}^{(1)}, \max_{t\in\mathfrak{R}(\delta)\cap [0,S]}X_{k}(t)\leq u_{k}^{(2)}, k=1,\cdots,p\right\}\\
&&= \sum_{k=1}^{p}P\left\{\max_{t\in [0,S]}X_{k}(t)>u_{k}^{(1)}\right\}+\sum_{k=1}^{p}P\left\{\max_{t\in\mathfrak{R}(\delta)\cap [0,S]}X_{k}(t)>u_{k}^{(2)}\right\}\\
&&\ \ \ -2\sum_{1\leq k<l\leq p}P\left\{\max_{t\in [0,S]}X_{k}(t)>u_{k}^{(1)}, \max_{t\in [0,S]}X_{l}(t)>u_{l}^{(1)}\right\}\\
&&\ \ \ -2\sum_{1\leq k<l\leq p}P\left\{\max_{t\in\mathfrak{R}(\delta)\cap [0,S]}X_{k}(t)>u_{k}^{(2)}, \max_{t\in \mathfrak{R}(\delta)\cap[0,S]}X_{l}(t)>u_{l}^{(2)}\right\}\\
&&\ \ \ -2\sum_{1\leq k<l\leq p}P\left\{\max_{t\in[0,S]}X_{k}(t)>u_{k}^{(1)}, \max_{t\in \mathfrak{R}(\delta)\cap[0,S]}X_{l}(t)>u_{l}^{(2)}\right\}\\
&&\ \ \ -\sum_{k=1}^{p}P\left\{\max_{t\in [0,S]}X_{k}(t)>u_{k}^{(1)}, \max_{t\in \mathfrak{R}(\delta)\cap[0,S]}X_{k}(t)>u_{k}^{(2)}\right\}+\sum_{3,p}\\
&&=:J_{1}+J_{2}-2J_{3}-2J_{4}-2J_{5}-J_{6}+\sum_{3,p}.
\end{eqnarray*}}
By Lemma 3 of Piterbarg (2004) and (\ref{eq3.4.2}), (\ref{eq3.4.3})
 as $T \to \IF$
\begin{eqnarray*}
\ccE{A_{1}}&\thicksim & \sum_{k=1}^{p}ST^{-1}(e^{-x_{k}-r_{kk}+\sqrt{2r_{kk}}z_{k}}+e^{-y_{k}-r_{kk}+\sqrt{2r_{kk}}z_{k}})\\
&=&\sum_{k=1}^{p}\frac{T^{a}}{T}(e^{-x_{k}-r_{kk}+\sqrt{2r_{kk}}z_{k}}+e^{-y_{k}-r_{kk}+\sqrt{2r_{kk}}z_{k}}).
\end{eqnarray*}
\ccE{With similar} arguments as for $A_{2}, A_{3}, A_{4}$ \cZ{in the
proof of Theorem 2.1}, we conclude that
$$\ccE{A_{k}=o(A_{1})},\ \ \ k=2,3,4.$$
\ccT{For the sum $A_{6}$, it is easy to see that each term in
$A_{6}$ can be bounded by $A_{3}$ or $A_{4}$, and thus $A_{6}=o(A_{1}).$} \he{The claim follows now easily borrowing the arguments in  p. 176 of Piterbarg (2004).}
\cE{\hfill $\Box$}

\textbf{Proof of Theorem 2.3}. First, recall that $\mathfrak{R}(\delta)$ is a dense grid in this case. By Lemma 5 of Piterbarg (2004) we have
\begin{eqnarray*}
&&\bigg|P\left\{a_{T}(M_{k}(T)-b_{T})\leq x_{k}, a_{T}(M_{k}^{\delta}(T)- b_{T})\leq y_{k}, k=1,\cdots,p\right\}\\
&&\ \ \ \ \ \ \ \ -P\left\{a_{T}(M_{k}(T)-b_{T})\leq x_{k}, a_{T}(M_{k}(T)- b_{T})\leq y_{k}, k=1,\cdots,p\right\}\bigg|\\
&&\leq\sum_{k=1}^{p}\bigg|P\left\{a_{T}(M_{k}^{\delta}(T)- b_{T})\leq y_{k}\right\}
      -P\left\{a_{T}(M_{k}(T)- b_{T})\leq y_{k}\right\}\bigg|\rightarrow0, \quad T\rightarrow\infty.
\end{eqnarray*}
 Since
\begin{eqnarray*}
&&P\left\{a_{T}(M_{k}(T)-b_{T})\leq x_{k}, a_{T}(M_{k}(T)- b_{T})\leq y_{k}, k=1,\cdots,p\right\}\\
&&=P\left\{a_{T}(M_{k}(T)-b_{T})\leq \min(x_{k}, y_{k}), k=1,\cdots,p\right\}
\end{eqnarray*}
in order to complete the proof, we only need to show that
\begin{eqnarray*}
&&\lim_{T\to \IF}P\left\{a_{T}(M_{k}(T)-b_{T})\leq \min(x_{k}, y_{k}), k=1,\cdots,p\right\}
=
 \int_{\mathbf{z} \in \mathbb{R}^p}
\cZ{\exp\left(-h(\mathbf{x},\mathbf{y},\mathbf{z})\right)}\Phidz,
\end{eqnarray*}
\ccE{which follows from Corollary 2.1}. \cE{\hfill $\Box$}

\section{Appendix}
\def\Dkl{\Delta_{kl}(nq,mq)}
\def\Dkk{\Delta_{kk}(nq,mq)}
\def\DklD{\Delta_{kl}(n\delta,m\delta)}
\def\DklQD{\Delta_{kl}(nq,m\delta)}
\def\DkkD{\Delta_{kk}(n\delta,m\delta)}
In this section, we give the detailed proof of Lemma 3.1 \iE{which is based on the results of \Quan{six} lemmas given below.}\\
\aa{Let in the following} $\mathcal{C}$ be constant whose value will
change from place to place. \aa{Define further}
\Quan{$r^{(h)}_{kl}(t,s)=hr_{kl}(t,s)+(1-h)\varrho_{kl}(t,s)$ for
$h\in[0,1]$ and $1\leq k,l\leq p$} and let
$\vartheta_{kk}(t)=\sup_{t<|nq-mq|\leq T}\{\varpi_{kk}(nq,mq)\}$,
where
$\varpi_{kk}(nq,mq)=\max\{r_{kk}(nq,mq),\varrho_{kk}(nq,mq)\}$.
Assumption (\ref{eq1.3}) \aa{implies that}
$\vartheta_{kk}(\epsilon)<1$ for \TT{all $T$} and any
$\epsilon\in(0,2^{-1/\alpha})$. Consequently, we may choose some
positive constant $\beta_{kk}$ such that \BQN \label{bkk}
\beta_{kk}<\frac{1-\vartheta_{kk}(\epsilon)}{1+\vartheta_{kk}(\epsilon)}
 \EQN
for all sufficiently large $T$. In the following we choose
$$0 < a < \aa{b}<\min_{k\in\{1,2,\cdots,p\}}\beta_{kk}$$
and we set $\Delta_{kl}(ns,mt):= |r_{kl}(ns,mt)-\varrho_{kl}(ns,mt)|$ for all possible indices $k,l$.

\textbf{Lemma 5.1}. {\sl Under conditions of Lemma 3.1, 
we have
$$\sum_{nq\in \mathcal{S}_{i},mq\in \mathcal{S}_{j}\atop i,j=1,2\cdots,n_{T}} \Dkk \int_{0}^{1}\frac{1}{\sqrt{1-r_{kk}^{(h)}(nq,mq)}}\exp\bigg(-\frac{(u_{k}^{(1)})^{2}}{1+r_{kk}^{(h)}(nq,mq)}\bigg)dh\rightarrow0$$
as $T\rightarrow\infty$.}

\textbf{Proof}:
 First, we consider the case t
$nq$ and $mq$ are in the same interval $\mathcal{S}$ which implies
$\varrho_{kk}(nq,mq)=r_{kk}(nq,mq)+(1-r_{kk}(nq,mq))\rho_{kk}(T)\sim
r_{kk}(nq,mq)$ as $T \to \IF$. Split the sum in the lemma into two parts as
\begin{eqnarray*}
\sum_{nq, mq\in \mathcal{S}_{i}\atop i=1,2\cdots,n_{T}, |nq-mq|<\epsilon}(\cdots)+\sum_{nq, mq\in \mathcal{S}_{i}\atop i=1,2\cdots,n_{T}, |nq-mq|\geq\epsilon}(\cdots)=:J_{T,1}+J_{T,2}.
\end{eqnarray*}
For the term $J_{T,1}$  note that Assumption (\ref{eq1.3}) implies for all
$|s-t|\leq\epsilon<2^{-1/\alpha}$
$$\frac{1}{2}|s-t|^{\alpha}\leq 1-r_{kk}(s,t)\leq 2|s-t|^{\alpha}.$$
\aa{By the  definition of $u_{k}^{(1)}$}
\begin{eqnarray}
\label{eq5.1.1}
(u_{k}^{(1)})^{2}=2\ln T-\ln\ln T+\frac{2}{\alpha}\ln\ln T+O(1), \quad T\to \IF.
\end{eqnarray}
Consequently, since further $q=\varepsilon(\ln T)^{-1/\alpha}$
\begin{eqnarray}
\label{eq5.1.2}
J_{T,1}&\leq&\mathcal{C}\sum_{nq, mq\in \mathcal{S}_{i}\atop i=1,2\cdots,n_{T}, |nq-mq|<\epsilon} \Dkk \frac{1}{\sqrt{1-r_{kk}(nq,mq)}}\exp\left(-\frac{(u_{k}^{(1)})^{2}}{1+r_{kk}(nq,mq)}\right)\nonumber\\
&\leq &\mathcal{C} \rho_{kk}(T)\sum_{nq, mq\in \mathcal{S}_{i}\atop i=1,2\cdots,n_{T}, |nq-mq|<\epsilon}
\frac{1-r_{kk}(nq,mq)}{\sqrt{1-r_{kk}(nq,mq)}}\exp\left(-\frac{(u_{k}^{(1)})^{2}}{2}\right)\exp\left(-\frac{(1-r_{kk}(nq,mq))(u_{k}^{(1)})^{2}}{2(1+r_{kk}(nq,mq))}\right)\nonumber\\
&\leq & \mathcal{C} \rho_{kk}(T)T^{-1}(\ln T)^{1/2-1/\alpha}\sum_{nq, mq\in \mathcal{S}_{i}\atop i=1,2\cdots,n_{T}, |nq-mq|<\epsilon}\sqrt{1-r_{kk}(nq,mq)}\exp\left(-\frac{(1-r_{kk}(nq,mq))(u_{k}^{(1)})^{2}}{2(1+r_{kk}(nq,mq))}\right)\nonumber\\
&\leq & \mathcal{C} \rho_{kk}(T)T^{-1}(\ln T)^{1/2-1/\alpha}\sum_{nq, mq\in \mathcal{S}_{i}\atop i=1,2\cdots,n_{T}, |nq-mq|<\epsilon}|nq-mq|^{\alpha/2}\exp\left(-\frac{1}{8}|nq-mq|^{\alpha}(u_{k}^{(1)})^{2}\right)\nonumber\\
&\leq & \mathcal{C} (\ln T)^{-1/2}\sum_{0<nq< \epsilon}(nq)^{\alpha/2}\exp\left(-\frac{1}{8}(nq)^{\alpha}(u_{k}^{(1)})^{2}\right)\nonumber\\
&\leq & \mathcal{C} (\ln T)^{-1/2}\sum_{0<nq< \epsilon}\exp\left(-\frac{1}{8}(nq)^{\alpha}\ln T\right)\nonumber\\
&\leq & \mathcal{C} (\ln T)^{-1/2}\sum_{n=1}^{\infty}\exp\left(-\frac{1}{8}(\varepsilon n)^{\alpha}\right)\nonumber\\
&\leq & \mathcal{C} (\ln T)^{-1/2}
\end{eqnarray}
 implying thus  $ \lim_{T \to \IF}J_{T,1}=0$. Using the fact that $u_{k}^{(1)}\thicksim
(2\log T)^{1/2}$ as $T\to \IF$ we obtain
\begin{eqnarray}
\label{eq5.1.3}
J_{T,2}&\leq&\mathcal{C}\sum_{nq, mq\in \mathcal{S}_{i}\atop i=1,2\cdots,n_{T}, |nq-mq|\geq\epsilon} \Dkk \exp\left(-\frac{(u_{k}^{(1)})^{2}}{1+\varpi_{kk}(nq,mq)}\right)\nonumber\\
&\leq & \mathcal{C}\sum_{nq, mq\in \mathcal{S}_{i}\atop i=1,2\cdots,n_{T}, |nq-mq|\geq\epsilon} \Dkk \exp\left(-\frac{(u_{k}^{(1)})^{2}}{1+\vartheta_{kk}(\epsilon)}\right)\nonumber\\
&\leq & \mathcal{C} \exp\left(-\frac{(u_{k}^{(1)})^{2}}{1+\vartheta_{kk}(\epsilon)}\right)\sum_{nq, mq\in \mathcal{S}_{i}\atop i=1,2\cdots,n_{T}, |nq-mq|\geq\epsilon}1\nonumber\\
&\leq & \mathcal{C} \frac{T}{q}T^{-\frac{2}{1+\vartheta_{kk}(\epsilon)}}\sum_{\epsilon<nq\leq T^{a}}1\nonumber\\
&\leq & \mathcal{C} T^{a-\frac{1-\vartheta_{_{kk}}(\epsilon)}{1+\vartheta_{kk}(\epsilon)}}(\ln T)^{2/\alpha}.
\end{eqnarray}
Since \Quan{$a<\min_{k\in\{1,2,\cdots,p\}}\beta_{kk}<\min_{k\in\{1,2,\cdots,p\}}\frac{1-\vartheta_{kk}(\epsilon)}{1+\vartheta_{kk}(\epsilon)}$} we have $\lim_{T\to \IF}J_{T,2}=0$. Second, we deal with the case that $nq\in \mathcal{S}_{i}$ and
$mq\in \mathcal{S}_{j}$, $i\neq j$. Note that in this case, the
distance between any two intervals $\mathcal{S}_{i}$ and
$\mathcal{S}_{j}$ is large than $T^{b}$. Split the sum into two
parts as
\begin{eqnarray*}
\sum_{nq\in\mathcal{S}_{i}, mq\in \mathcal{S}_{j},i\neq j\atop i,j=1,2\cdots,n_{T}, |nq-mq|<T^{\beta_{kk}}}(\cdots)
+\sum_{nq\in\mathcal{S}_{i}, mq\in \mathcal{S}_{j}i\neq j\atop i,j=1,2\cdots,n_{T}, |nq-mq|\geq T^{\beta_{kk}}}(\cdots)=:J_{T,3}+J_{T,4}.
\end{eqnarray*}
Similarly to the derivation of (\ref{eq5.1.3}), for large enough $T$
we have
\begin{eqnarray}
\label{eq5.1.4}
J_{T,3}&\leq&\mathcal{C}\sum_{nq\in\mathcal{S}_{i}, mq\in \mathcal{S}_{j}i\neq j\atop i,j=1,2\cdots,n_{T}, |nq-mq|<T^{\beta_{kk}}} \Dkk \exp\left(-\frac{(u_{k}^{(1)})^{2}}{1+\varpi_{kk}(nq,mq)}\right)\nonumber\\
&\leq & \mathcal{C}\sum_{nq\in\mathcal{S}_{i}, mq\in \mathcal{S}_{j}i\neq j\atop i,j=1,2\cdots,n_{T}, |nq-mq|<T^{\beta_{kk}}} \Dkk \exp\left(-\frac{(u_{k}^{(1)})^{2}}{1+\vartheta_{kk}(T^{b})}\right)\nonumber\\
&\leq & \mathcal{C} \exp\left(-\frac{(u_{k}^{(1)})^{2}}{1+\vartheta_{kk}(\epsilon)}\right)\sum_{nq\in\mathcal{S}_{i}, mq\in \mathcal{S}_{j}i\neq j\atop i,j=1,2\cdots,n_{T}, |nq-mq|<T^{\beta_{kk}}}1\nonumber\\
&\leq & \mathcal{C} \frac{T}{q}T^{-\frac{2}{1+\vartheta_{kk}(\epsilon)}}\sum_{T^{b}<nq\leq T^{\beta_{kk}}}1\nonumber\\
&\leq & \mathcal{C}T^{\beta_{kk}-\frac{1-\vartheta_{kk}(\epsilon)}{1+\vartheta_{kk}(\epsilon)}}(\ln T)^{2/\alpha}.
\end{eqnarray}
Consequently,  \aa{by \eqref{bkk}} $J_{T,3}\rightarrow0$ as $T\rightarrow\infty$. By the Assumption (\ref{eq1.4}) we have $\vartheta_{kk}(t)\ln t\leq
K$ for all sufficiently large $t$ and some constant $K$. Thus,
$\varpi_{kk}(nq,mq)\leq \vartheta_{kk}(T^{\beta_{kk}})\leq K/\ln
T^{\beta_{kk}}$ for $|nq-mq|>T^{\beta_{kk}}$. Now using
(\ref{eq5.1.1}) again for all large $T$ we obtain
\begin{eqnarray}
\label{eq5.1.5}
\frac{T^{2}}{q^{2}\ln T}\exp\left(-\frac{(u_{k}^{(1)})^{2}}{1+\vartheta_{kk}(T^{\beta_{kk}})}\right)
&\leq& \frac{T^{2}}{q^{2}\ln T}\exp\left(-\frac{(u_{k}^{(1)})^{2}}{1+K/\ln T^{\beta_{kk}}}\right)\nonumber\\
&=&\frac{T^{2}}{q^{2}\ln T}\left(T^{-2}\ln T (\ln T)^{-2/\alpha}\right)^{\frac{1}{1+K/\ln T^{\beta_{kk}}}} (1+ o(1))\nonumber\\
&=&\varepsilon^{-2}T^{(2K/\ln T^{\beta_{kk}})/(1+K/\ln T^{\beta_{k}})}(\ln T)^{((2/\alpha_{kk}-1)K/\ln T^{\beta_{kk}})/(1+K/\ln T^{\beta_{kk}})}\nonumber\\
&=&O(1).
\end{eqnarray}
Further, we have
\begin{eqnarray}
\label{eq5.1.6}
J_{T,4}&\leq&\mathcal{C}\sum_{nq\in\mathcal{S}_{i}, mq\in \mathcal{S}_{j}i\neq j\atop i,j=1,2\cdots,n_{T}, |nq-mq|\geq T^{\beta_{kk}}} \Dkk \exp\left(-\frac{(u_{k}^{(1)})^{2}}{1+\varpi_{kk}(nq,mq)}\right)\nonumber\\
&\leq&\mathcal{C}\exp\left(-\frac{(u_{k}^{(1)})^{2}}{1+\vartheta_{kk}(T^{\beta_{kk}})}\right)\sum_{nq\in\mathcal{S}_{i}, mq\in \mathcal{S}_{j}i\neq j\atop i,j=1,2\cdots,n_{T}, |nq-mq|\geq T^{\beta_{kk}}} \Dkk \nonumber\\
&=& \mathcal{C}\frac{T^{2}}{q^{2}\ln T}\exp\left(-\frac{(u_{k}^{(1)})^{2}}{1+\vartheta_{kk}(T^{\beta_{kk}})}\right)\cdot \frac{q^{2}\ln T}{T^{2}}\sum_{nq\in\mathcal{S}_{i}, mq\in \mathcal{S}_{j}i\neq j\atop i,j=1,2\cdots,n_{T}, |nq-mq|\geq T^{\beta_{kk}}} \Dkk \nonumber\\
&\leq& \mathcal{C} \frac{q^{2}\ln T}{T^{2}}\sum_{nq\in\mathcal{S}_{i}, mq\in \mathcal{S}_{j}i\neq j\atop i,j=1,2\cdots,n_{T}, |nq-mq|\geq T^{\beta_{kk}}}\left|r_{kk}(nq,mq)-\frac{r_{kk}}{\ln T}\right|\nonumber\\
&\leq&  \mathcal{C}\frac{q^{2}}{\beta^{kk}T^{2}}\sum_{nq\in\mathcal{S}_{i}, mq\in \mathcal{S}_{j}i\neq j\atop i,j=1,2\cdots,n_{T}, |nq-mq|\geq T^{\beta_{kk}}}|r_{kk}(nq,mq)\ln(|nq-mq|)-r_{kk}|\nonumber\\
&&+\mathcal{C}r_{kk}\frac{q^{2}}{T^{2}}\sum_{nq\in\mathcal{S}_{i}, mq\in \mathcal{S}_{j}i\neq j\atop i,j=1,2\cdots,n_{T}, |nq-mq|\geq T^{\beta_{kk}}}\left|1-\frac{\ln T}{\ln(|nq-mq|)}\right|\nonumber\\
&\leq&  \mathcal{C}\frac{q}{\beta^{kk}T}\sum_{T^{\beta_{kk}} \leq nq\leq T }|r_{kk}(nq)\ln(nq)-r_{kk}|
+\mathcal{C}r_{kk}\frac{q}{T}\sum_{T^{\beta_{kk}} \leq nq\leq T }\left|1-\frac{\ln T}{\ln(nq)}\right|.
\end{eqnarray}
By the Assumption (\ref{eq1.4}) the first term on the
right-hand-side of (\ref{eq5.1.6}) tends to 0 as $T\to \infty$.
Furthermore, the second term of the right-hand-side of
(\ref{eq5.1.6}) also tends to 0 by an integral estimate as in the
proof of Lemma 6.4.1 of Leadbetter et al.\ (1983). Now from
(\ref{eq5.1.2})-(\ref{eq5.1.6}), we get that the sum in the claim of the lemma tends to $0$ as $T\rightarrow\infty$.\cE{\hfill $\Box$}

\textbf{Lemma 5.2}. {\sl Under conditions of Lemma 3.1, we have
$$\sum_{n\delta\in \mathcal{S}_{i},m\delta\in \mathcal{S}_{j}\atop i,j=1,2\cdots,n_{T}} \DkkD \int_{0}^{1}\frac{1}{\sqrt{1-r_{kk}^{(h)}(n\delta,m\delta)}}\exp\bigg(-\frac{(u_{k}^{(2)})^{2}}{1+r_{kk}^{(h)}(n\delta,m\delta)}\bigg)dh\rightarrow0$$
as $T\rightarrow\infty$.}

\textbf{Proof:} The claim is established by following very closely the proof of Lemma 5.1.  \cE{\hfill $\Box$}

\textbf{Lemma 5.3}. {\sl Under conditions of Lemma 3.1, we have
$$\sum_{nq\in \mathcal{S}_{i},m\delta\in \mathcal{S}_{j}, nq\neq m\delta\atop i,j=1,2\cdots,n_{T}} \DklQD  \int_{0}^{1}\frac{1}{\sqrt{1-r_{kk}^{(h)}(nq,m\delta)}}\exp\bigg(-\frac{(u_{k}^{(1)})^{2}+(u_{k}^{(2)})^{2}}{2(1+r_{kk}^{(h)}(nq,m\delta))}\bigg)dh\rightarrow0$$
as $T\rightarrow\infty$.}

\textbf{Proof:} We only prove the case that \aa{$\mathfrak{R}(\delta)$} is a sparse
grid, since the proof of the Pickands grid is the same.

First, we consider the case that $nq,m\delta$ in the same
interval $\mathcal{S}$. Note that in this case,
$\varrho_{kk}(nq,m\delta)=r_{kk}(nq,m\delta)+(1-r_{kk}(nq,m\delta))\rho_{kk}(T)\sim
r_{kk}(nq,m\delta)$ for sufficiently large $T$. Split the sum into
two parts as
\begin{eqnarray*}
\sum_{nq, m\delta\in \mathcal{S}_{i}\atop i=1,2\cdots,n_{T}, |nq-m\delta|<\epsilon} (\cdots)+\sum_{nq, m\delta\in \mathcal{S}_{i}\atop i=1,2\cdots,n_{T}, |nq-m\delta|\geq\epsilon}(\cdots)=:S_{T,1}+S_{T,2}.
\end{eqnarray*}
The definition of $u_{k}^{(2)}$ implies thus
\begin{eqnarray}
\label{eq5.3.1}
(u_{k}^{(2)})^{2}=2\ln T-\ln\ln T+2\ln \delta^{-1}+O(1)
\end{eqnarray}
\def\ukAB{u_{k,12}}
for sparse grids. Consequently, since further $q=\varepsilon(\ln
T)^{-1/\alpha}$ and the definition of $\delta$ we obtain (set $\ukAB:= (u_{k}^{(1)})^{2}+(u_{k}^{(2)})^2$)
\begin{eqnarray}
\label{eq5.3.2}
S_{T,1}&\leq&\mathcal{C}\sum_{nq, m\delta\in \mathcal{S}_{i}\atop i=1,2\cdots,n_{T}, |nq-m\delta|<\epsilon} \DklQD  \frac{1}{\sqrt{1-r_{kk}(nq,m\delta)}}
\exp\left(-\frac{\ukAB}{2(1+r_{kk}(nq,m\delta)}\right)\nonumber\\
&\leq &\mathcal{C} \rho_{kk}(T)\sum_{nq, m\delta\in \mathcal{S}_{i}\atop i=1,2\cdots,n_{T}, |nq-m\delta|<\epsilon}
\frac{1-r_{kk}(nq,m\delta)}{\sqrt{1-r_{kk}(nq,m\delta)}}\exp\left(-\frac{\ukAB}{4}\right)\exp\left(-\frac{(1-r_{kk}(nq,m\delta))
\ukAB}{4(1+r_{kk}(nq,m\delta))}\right)\nonumber\\
&\leq & \mathcal{C} \rho_{kk}(T)T^{-1}(\ln T)^{1/2-1/2\alpha}\delta^{1/2}\sum_{nq, m\delta\in \mathcal{S}_{i}\atop i=1,2\cdots,n_{T}, |nq-m\delta|<\epsilon}\sqrt{1-r_{kk}(nq,m\delta)}\exp\left(-\frac{(1-r_{kk}(nq,m\delta))\ukAB}{4(1+r_{kk}(nq,m\delta))}\right)\nonumber\\
&\leq & \mathcal{C} \rho_{kk}(T)T^{-1}(\ln T)^{1/2-1/2\alpha}\delta^{1/2}\sum_{nq, m\delta\in \mathcal{S}_{i}\atop i=1,2\cdots,n_{T}, |nq-m\delta|<\epsilon}|nq-m\delta|^{\alpha/2}
\exp\left(-\frac{1}{16}|nq-m\delta|^{\alpha}\ukAB\right)\nonumber\\
&\leq & \mathcal{C} T^{-1}(\ln T)^{-1/2-1/2\alpha}\delta^{1/2}\sum_{nq, m\delta\in \mathcal{S}_{i}\atop i=1,2\cdots,n_{T}, |nq-m\delta|<\epsilon}
\exp\left(-\frac{1}{16}|nq-m\delta|^{\alpha}\ln T\right).
\end{eqnarray}
Since \aa{$\mathfrak{R}(\delta)$} is a sparse grid, $\delta(\ln
T)^{1/\alpha}\rightarrow\infty$. A simple calculation shows that
\begin{eqnarray*}
\sum_{nq, m\delta\in
\mathcal{S}_{i}\atop i=1,2\cdots,n_{T}, |nq-m\delta|<\epsilon}
\exp\left(-\frac{1}{16}|nq-m\delta|^{\alpha}\ln T\right)\leq \mathcal{C}T\delta^{-1}\sum_{0<nq<\epsilon}
\exp\left(-\frac{1}{16}(nq)^{\alpha}\ln T\right)\leq \mathcal{C}T\delta^{-1}.
\end{eqnarray*}
Hence the assumption $\delta(\ln T)^{1/\alpha}\rightarrow\infty$ implies  $S_{T,1}\leq \mathcal{C}(\ln
T)^{-1/2-1/2\alpha}\delta^{-1/2} = o( (\ln T)^{-1/2})$. \aa{Since  $u_{k}^{(i)}\thicksim (2\log T)^{1/2},i=1,2$ }
\begin{eqnarray}
\label{eq5.3.3}
S_{T,2}&\leq&\mathcal{C}\sum_{nq, m\delta\in \mathcal{S}_{i}\atop i=1,2\cdots,n_{T}, |nq-m\delta|\geq\epsilon} \DklQD  \exp\left(-\frac{\ukAB}{2(1+\varpi_{kk}(nq,m\delta))}\right)\nonumber\\
&\leq & \mathcal{C}\sum_{nq, m\delta\in \mathcal{S}_{i}\atop i=1,2\cdots,n_{T}, |nq-m\delta|\geq\epsilon} \DklQD  \exp\left(-\frac{\ukAB}{2(1+\vartheta_{kk}(\epsilon))}\right)\nonumber\\
&\leq & \mathcal{C} \exp\left(-\frac{\ukAB}{2(1+\vartheta_{kk}(\epsilon))}\right)\sum_{nq, m\delta\in \mathcal{S}_{i}\atop i=1,2\cdots,n_{T}, |nq-m\delta|\geq\epsilon}1\nonumber\\
&\leq & \mathcal{C} \frac{T}{q}T^{-\frac{2}{1+\vartheta_{kk}(\epsilon)}}\sum_{0<m\delta\leq T^{a}}1\nonumber\\
&\leq & \mathcal{C} T^{a-\frac{1-\vartheta_{_{kk}}(\epsilon)}{1+\vartheta_{kk}(\epsilon)}}\delta^{-1}(\ln T)^{1/\alpha}.
\end{eqnarray}
\aa{In view of \eqref{bkk} and  \TT{the fact that
$\mathfrak{R}(\delta)$ is a sparse grid},
$ \lim_{T \to \IF} S_{T,2}=0$.}\\
Second, we deal with the case that $nq\in \mathcal{S}_{i}$ and
$mq\in \mathcal{S}_{j}$, $i\neq j$. Again we split the sum into two parts as
\begin{eqnarray*}
\sum_{nq\in\mathcal{S}_{i}, m\delta\in \mathcal{S}_{j},i\neq j\atop i,j=1,2\cdots,n_{T}, |nq-m\delta|<T^{\beta_{kk}}} (\cdots)
+\sum_{nq\in\mathcal{S}_{i}, m\delta\in \mathcal{S}_{j}i\neq j\atop i,j=1,2\cdots,n_{T}, |nq-m\delta|\geq T^{\beta_{kk}}}(\cdots)=:S_{T,3}+S_{T,4}.
\end{eqnarray*}
Similarly to the derivation of (\ref{eq5.3.3}) for large enough $T$
we have
\begin{eqnarray}
\label{eq5.3.4}
S_{T,3}&\leq&\mathcal{C}\sum_{nq\in\mathcal{S}_{i}, m\delta\in \mathcal{S}_{j}i\neq j\atop i,j=1,2\cdots,n_{T}, |nq-m\delta|<T^{\beta_{kk}}} \DklQD  \exp\left(-\frac{\ukAB}{2(1+\varpi_{kk}(nq,m\delta))}\right)\nonumber\\
&\leq & \mathcal{C}\sum_{nq\in\mathcal{S}_{i}, m\delta\in \mathcal{S}_{j}i\neq j\atop i,j=1,2\cdots,n_{T}, |nq-m\delta|<T^{\beta_{kk}}} \DklQD  \exp\left(-\frac{\ukAB}{2(1+\vartheta_{kk}(T^{b}))}\right)\nonumber\\
&\leq & \mathcal{C} \exp\left(-\frac{\ukAB}{2(1+\vartheta_{kk}(\epsilon))}\right)\sum_{nq\in\mathcal{S}_{i}, m\delta\in \mathcal{S}_{j}i\neq j\atop i,j=1,2\cdots,n_{T}, |nq-m\delta|<T^{\beta_{kk}}}1\nonumber\\
&\leq & \mathcal{C} \frac{T}{q}T^{-\frac{2}{1+\vartheta_{kk}(\epsilon)}}\sum_{0<n\delta\leq T^{\beta_{kk}}}1\nonumber\\
&\leq & \mathcal{C}T^{\beta_{kk}-\frac{1-\vartheta_{kk}(\epsilon)}{1+\vartheta_{kk}(\epsilon)}}\aa{\delta^{-1}}(\ln T)^{1/\alpha}.
\end{eqnarray}
Consequently, by \eqref{bkk} $\lim_{T\to \IF} S_{T,3}=0$. By the Assumption (\ref{eq1.4}) we have $\vartheta_{kk}(t)\ln t\leq
K$ for all sufficiently large $t$ and some constant $K$. Thus,
$\varpi_{kk}(nq,m\delta)\leq \vartheta_{kk}(T^{\beta_{kk}})\leq
K/\ln T^{\beta_{kk}}$ for $|nq-m\delta|>T^{\beta_{kk}}$. Now using
(\ref{eq5.1.1}) and (\ref{eq5.3.1}) again for all large $T$ and
$|nq-m\delta|>T^{\beta_{kk}}$ we obtain
\begin{eqnarray}
\label{eq5.3.5}
\frac{T^{2}}{q\delta\ln T}\exp\left(-\frac{\ukAB}{2(1+\vartheta_{kk}(T^{\beta_{kk}}))}\right)
&\leq& \frac{T^{2}}{q\delta\ln T}\exp\left(-\frac{\ukAB}{2(1+K/\ln T^{\beta_{kk}})}\right)\nonumber\\
&=&\frac{T^{2}}{q\delta\ln T}\left(T^{-2}\ln T (\ln T)^{-1/\alpha}\delta\right)^{\frac{1}{1+K/\ln T^{\beta_{kk}}}} (1+ o(1))\nonumber\\
&=&O(1).
\end{eqnarray}
Now, with similar arguments as in the proof of (\ref{eq5.1.6}) we obtain
\begin{eqnarray}
\label{eq5.3.6}
S_{T,4}&\leq&\mathcal{C}\sum_{nq\in\mathcal{S}_{i}, m\delta\in \mathcal{S}_{j}i\neq j\atop i,j=1,2\cdots,n_{T}, |nq-m\delta|\geq T^{\beta_{kk}}} \DklQD  \exp\left(-\frac{\ukAB}{2(1+\varpi_{kk}(nq,m\delta))}\right)\nonumber\\
&\leq&\mathcal{C}\exp\left(-\frac{\ukAB}{2(1+\vartheta_{kk}(T^{\beta_{kk}}))}\right)\sum_{nq\in\mathcal{S}_{i}, m\delta\in \mathcal{S}_{j}i\neq j\atop i,j=1,2\cdots,n_{T}, |nq-m\delta|\geq T^{\beta_{kk}}} \DklQD  \nonumber\\
&=& \mathcal{C}\frac{T^{2}}{q\delta\ln T}\exp\left(-\frac{\ukAB}{2(1+\vartheta_{kk}(T^{\beta_{kk}}))}\right) \frac{q\delta\ln T}{T^{2}}\sum_{nq\in\mathcal{S}_{i}, m\delta\in \mathcal{S}_{j}i\neq j\atop i,j=1,2\cdots,n_{T}, |nq-m\delta|\geq T^{\beta_{kk}}} \DklQD  \nonumber\\
&\leq& \mathcal{C} \frac{q\delta\ln T}{T^{2}}\sum_{nq\in\mathcal{S}_{i}, m\delta\in \mathcal{S}_{j}i\neq j\atop i,j=1,2\cdots,n_{T}, |nq-m\delta|\geq T^{\beta_{kk}}}\left|r_{kk}(nq,m\delta)-\frac{r_{kk}}{\ln T}\right|\nonumber\\
&\leq&  \mathcal{C}\frac{q\delta}{\beta^{kk}T^{2}}\sum_{nq\in\mathcal{S}_{i}, m\delta\in \mathcal{S}_{j}i\neq j\atop i,j=1,2\cdots,n_{T}, |nq-m\delta|\geq T^{\beta_{kk}}}|r_{kk}(nq,m\delta)\ln(|nq-m\delta|)-r_{kk}|\nonumber\\
&&+\mathcal{C}r_{kk}\frac{q\delta}{T^{2}}\sum_{nq\in\mathcal{S}_{i}, m\delta\in \mathcal{S}_{j}i\neq j\atop i,j=1,2\cdots,n_{T}, |nq-m\delta|\geq T^{\beta_{kk}}}\left|1-\frac{\ln T}{\ln(|nq-m\delta|)}\right|.
\end{eqnarray}
By the Assumption (\ref{eq1.4}) the first term on the
right-hand-side of (\ref{eq5.3.6}) tends to 0 as $T\to \infty$.
Furthermore, the second term of the right-hand-side of
(\ref{eq5.3.6}) also tends to 0 by an integral estimate as in the
proof of Lemma 6.4.1 of Leadbetter et al.\ (1983). Now the claim follows from
(\ref{eq5.3.2})-(\ref{eq5.3.6}).\cE{\hfill
$\Box$}

\textbf{Lemma 5.4}. {\sl Under conditions of Lemma 3.1, we have for
$k<l$
$$\sum_{nq\in \mathcal{S}_{i},mq\in \mathcal{S}_{j}\atop i,j=1,2\cdots,n_{T}} \Dkl\int_{0}^{1}\frac{1}{\sqrt{1-r_{kl}^{(h)}(nq,mq)}}\exp\bigg(-\frac{(u_{k}^{(1)})^{2}+(u_{l}^{(1)})^{2}}{2(1+r_{kl}^{(h)}(nq,mq))}\bigg)dh\rightarrow0$$
as $T\rightarrow\infty$.}

\textbf{Proof:}  Let $\vartheta_{kl}(t)=\sup_{|nq-mq|\geq
t}\{\varpi_{kl}(nq,mq)\}$, where
$\varpi_{kl}(nq,mq)=\max\{r_{kl}(nq,mq),\varrho_{kl}(nq,mq)\}$. From
Assumption (\ref{eq1.5}) and the definition of
$\varrho_{kl}(nq,mq)$, we have $\vartheta_{kl}(0)<1$ for \TT{all
$T$}. Consequently, we may choose some positive constant
$\beta_{kl}$ such that
$\beta_{kl}<\frac{1-\vartheta_{kl}(0)}{1+\vartheta_{kl}(0)}$ for all
sufficiently large $T$.  Split the sum into two parts as
\begin{eqnarray*}
\sum_{nq\in\mathcal{S}_{i}, mq\in \mathcal{S}_{j}\atop i,j=1,2\cdots,n_{T}, |nq-mq|<T^{\beta_{kl}}}(\cdots)
+\sum_{nq\in\mathcal{S}_{i}, mq\in \mathcal{S}_{j}\atop i,j=1,2\cdots,n_{T}, |nq-mq|\geq T^{\beta_{kl}}}(\cdots)=:R_{T,1}+R_{T,2}.
\end{eqnarray*}
Similarly to the derivation of (\ref{eq5.1.3}), for large enough $T$ we have
\begin{eqnarray*}
\label{eq3.1.22}
R_{T,1}&\leq&\mathcal{C}\sum_{nq\in\mathcal{S}_{i}, mq\in \mathcal{S}_{j}\atop i,j=1,2\cdots,n_{T}, |nq-mq|<T^{\beta_{kl}}} \Dkl\exp\left(-\frac{(u_{k}^{(1)})^{2}+(u_{l}^{(1)})^{2}}{2(1+\varpi_{kl}(nq,mq))}\right)\nonumber\\
&\leq & \mathcal{C}\sum_{nq\in\mathcal{S}_{i}, mq\in \mathcal{S}_{j}\atop i,j=1,2\cdots,n_{T}, |nq-mq|<T^{\beta_{kl}}} \Dkl\exp\left(-\frac{(u_{k}^{(1)})^{2}+(u_{l}^{(1)})^{2}}{2(1+\vartheta_{kl}(0))}\right)\nonumber\\
&\leq & \mathcal{C} \exp\left(-\frac{(u_{k}^{(1)})^{2}+(u_{l}^{(1)})^{2}}{2(1+\vartheta_{kl}(0))}\right)\sum_{nq\in\mathcal{S}_{i}, mq\in \mathcal{S}_{j}\atop i,j=1,2\cdots,n_{T}, |nq-mq|<T^{\beta_{kl}}}1\nonumber\\
&\leq & \mathcal{C} \frac{T}{q}T^{-\frac{2}{1+\vartheta_{kl}(0)}}\sum_{0<nq\leq T^{\beta_{kl}}}1\nonumber\\
&\leq & \mathcal{C}T^{\beta_{kl}-\frac{1-\vartheta_{kl}(0)}{1+\vartheta_{kl}(0)}}(\ln T)^{2/\alpha}.
\end{eqnarray*}
Consequently,  $R_{T,1}\rightarrow0$ as $T\rightarrow\infty$ which follows by the fact that  $\beta_{kl}<\frac{1-\vartheta_{kl}(0)}{1+\vartheta_{kl}(0)}$.\\
By the Assumption (\ref{eq1.4}) we have $\vartheta_{kl}(t)\ln t\leq
K$ for all sufficiently large $t$. Thus, $\varpi_{kl}(nq,mq)\leq
\vartheta_{kl}(T^{\beta_{kl}})\leq K/\ln T^{\beta_{kl}}$ for
$|nq-mq|>T^{\beta_{kl}}$. Now with the similar arguments as for
(\ref{eq5.1.5}) we obtain
\begin{eqnarray*}
\frac{T^{2}}{q^{2}\ln T}\exp\left(-\frac{(u_{k}^{(1)})^{2}+(u_{l}^{(1)})^{2}}{2(1+\vartheta_{kl}(T^{\beta_{kl}}))}\right)=O(1).
\end{eqnarray*}
Thus, for $R_{T,2}$ we have
\begin{eqnarray*}
\label{eq3.1.3}
R_{T,2}&\leq&\mathcal{C}\sum_{nq\in\mathcal{S}_{i}, mq\in \mathcal{S}_{j}\atop i,j=1,2\cdots,n_{T}, |nq-mq|\geq T^{\beta_{kl}}} \Dkl\exp\left(-\frac{(u_{k}^{(1)})^{2}+(u_{l}^{(1)})^{2}}{2(1+\varpi_{kl}(nq,mq))}\right)\nonumber\\
&\leq&\mathcal{C}\exp\left(-\frac{(u_{k}^{(1)})^{2}+(u_{l}^{(1)})^{2}}{2(1+\vartheta_{kl}(T^{\beta_{kl}}))}\right)\sum_{nq\in\mathcal{S}_{i}, mq\in \mathcal{S}_{j}\atop i,j=1,2\cdots,n_{T}, |nq-mq|\geq T^{\beta_{kl}}} \Dkl\nonumber\\
&\leq& \mathcal{C}\frac{T^{2}}{q^{2}\ln T}\exp\left(-\frac{(u_{k}^{(1)})^{2}+(u_{l}^{(1)})^{2}}{2(1+\vartheta_{kl}(T^{\beta_{kl}}))}\right)\cdot \frac{q^{2}\ln T}{T^{2}}\sum_{nq\in\mathcal{S}_{i}, mq\in \mathcal{S}_{j}\atop i,j=1,2\cdots,n_{T}, |nq-mq|\geq T^{\beta_{kl}}} \Dkl\nonumber\\
&\leq& \mathcal{C} \frac{q^{2}\ln T}{T^{2}}\sum_{nq\in\mathcal{S}_{i}, mq\in \mathcal{S}_{j}\atop i,j=1,2\cdots,n_{T}, |nq-mq|\geq T^{\beta_{kl}}}\left|r_{kl}(nq,mq)-\frac{r_{kl}}{\ln T}\right|.
\end{eqnarray*}
By the same arguments as those in Lemma 5.1, we have $\lim_{T\to \IF} R_{T,2}=0$ and thus the claim follows. \cE{\hfill $\Box$}

\textbf{Lemma 5.5}. {\sl Under conditions of Lemma 3.1, we have for
$k<l$
$$\sum_{n\delta\in \mathcal{S}_{i},m\delta\in \mathcal{S}_{j}\atop i,j=1,2\cdots,n_{T}} \DklD \int_{0}^{1}\frac{1}{\sqrt{1-r_{kl}^{(h)}(n\delta,m\delta)}}\exp\bigg(-\frac{(u_{k}^{(2)})^{2}+(u_{l}^{(2)})^{2}}{2(1+r_{kl}^{(h)}(n\delta,m\delta))}\bigg)dh\rightarrow0$$
as $T\rightarrow\infty$.}

\textbf{Proof:} The claim follows with the same arguments as in the proof of Lemma 5.3. \hfill $\Box$

\textbf{Lemma 5.6}. {\sl Under conditions of Lemma 3.1, we have for $k<l$
$$\sum_{nq\in \mathcal{S}_{i},m\delta\in \mathcal{S}_{j},nq\neq m\delta\atop i,j=1,2\cdots,n_{T}} \DklQD  \int_{0}^{1}\frac{1}{\sqrt{1-r_{kl}^{(h)}(nq,m\delta)}}\exp\bigg(-\frac{(u_{k}^{(1)})^{2}+(u_{l}^{(2)})^{2}}{2(1+r_{kl}^{(h)}(nq,m\delta))}\bigg)dh\rightarrow0$$
as $T\rightarrow\infty$.}

\textbf{Proof:}  As in Lemma 5.2, we also only prove the case that
\aa{$\mathfrak{R}(\delta)$} is a sparse grid. Split the sum into two parts as
\Quan{(with the same definition $T^{\beta_{kl}}$ as in Lemma 5.3.)}
\begin{eqnarray*}
\label{eq3.1.11}
\sum_{nq\in\mathcal{S}_{i}, m\delta\in \mathcal{S}_{j}\atop i,j=1,2\cdots,n_{T}, |nq-m\delta|<T^{\beta_{kl}}}(\cdots)
+\sum_{nq\in\mathcal{S}_{i}, m\delta\in \mathcal{S}_{j}\atop i,j=1,2\cdots,n_{T}, |nq-m\delta|\geq T^{\beta_{kl}}}(\cdots)=:M_{T,1}+M_{T,2}.
\end{eqnarray*}
Similarly to the derivation of (\ref{eq5.1.3}), for large enough $T$
we have
\begin{eqnarray*}
\label{eq3.1.22}
M_{T,1}&\leq&\mathcal{C}\sum_{nq\in\mathcal{S}_{i}, m\delta\in \mathcal{S}_{j}\atop i,j=1,2\cdots,n_{T}, |nq-m\delta|<T^{\beta_{kl}}} \DklQD  \exp\left(-\frac{(u_{k}^{(1)})^{2}+(u_{l}^{(2)})^{2}}{2(1+\varpi_{kl}(nq,m\delta))}\right)\nonumber\\
&\leq & \mathcal{C}\sum_{nq\in\mathcal{S}_{i}, m\delta\in \mathcal{S}_{j}\atop i,j=1,2\cdots,n_{T}, |nq-m\delta|<T^{\beta_{kl}}} \DklQD  \exp\left(-\frac{(u_{k}^{(1)})^{2}+(u_{l}^{(2)})^{2}}{2(1+\vartheta_{kl}(0))}\right)\nonumber\\
&\leq & \mathcal{C} \exp\left(-\frac{(u_{k}^{(1)})^{2}+(u_{l}^{(2)})^{2}}{2(1+\vartheta_{kl}(0))}\right)\sum_{nq\in\mathcal{S}_{i}, m\delta\in \mathcal{S}_{j}\atop i,j=1,2\cdots,n_{T}, |nq-m\delta|<T^{\beta_{kl}}}1\nonumber\\
&\leq & \mathcal{C} \frac{T}{q}T^{-\frac{2}{1+\vartheta_{kl}(0)}}\sum_{0<m\delta\leq T^{\beta_{kl}}}1\nonumber\\
&\leq & \mathcal{C}T^{\beta_{kl}-\frac{1-\vartheta_{kl}(0)}{1+\vartheta_{kl}(0)}}(\ln T)^{1/\alpha}\delta^{-1}.
\end{eqnarray*}
Thus, $M_{T,1}\rightarrow0$ as $T\rightarrow\infty$ from the facts
that  $\beta_{kl}<\frac{1-\vartheta_{kl}(0)}{1+\vartheta_{kl}(0)}$
\Quan{and $\delta(\ln
T)^{1/\alpha}\rightarrow\infty$}. As for (\ref{eq5.3.5}) we have
\begin{eqnarray*}
\frac{T^{2}}{q\delta\ln T}\exp\left(-\frac{(u_{k}^{(1)})^{2}+(u_{l}^{(2)})^{2}}{2(1+\vartheta_{kl}(T^{\beta_{kl}}))}\right)=O(1).
\end{eqnarray*}
Consequently
\begin{eqnarray*}
\label{eq3.1.3}
M_{T,2}&\leq&\mathcal{C}\sum_{nq\in\mathcal{S}_{i}, m\delta\in \mathcal{S}_{j}\atop i,j=1,2\cdots,n_{T}, |nq-m\delta|\geq T^{\beta_{kl}}} \DklQD  \exp\left(-\frac{(u_{k}^{(1)})^{2}+(u_{l}^{(2)})^{2}}{2(1+\varpi_{kl}(nq,m\delta))}\right)\nonumber\\
&\leq&\mathcal{C}\exp\left(-\frac{(u_{k}^{(1)})^{2}+(u_{l}^{(2)})^{2}}{2(1+\vartheta_{kl}(T^{\beta_{kl}}))}\right)\sum_{nq\in\mathcal{S}_{i}, m\delta\in \mathcal{S}_{j}\atop i,j=1,2\cdots,n_{T}, |nq-m\delta|\geq T^{\beta_{kl}}} \DklQD  \nonumber\\
&\leq& \mathcal{C}\frac{T^{2}}{q\delta\ln T}\exp\left(-\frac{(u_{k}^{(1)})^{2}+(u_{l}^{(2)})^{2}}{2(1+\vartheta_{kl}(T^{\beta_{kl}}))}\right)\cdot \frac{q\delta\ln T}{T^{2}}\sum_{nq\in\mathcal{S}_{i}, m\delta\in \mathcal{S}_{j}\atop i,j=1,2\cdots,n_{T}, |nq-m\delta|\geq T^{\beta_{kl}}} \DklQD  \nonumber\\
&\leq& \mathcal{C} \frac{q\delta\ln T}{T^{2}}\sum_{nq\in\mathcal{S}_{i}, m\delta\in \mathcal{S}_{j}\atop i,j=1,2\cdots,n_{T}, |nq-m\delta|\geq T^{\beta_{kl}}}\left|r_{kl}(nq,m\delta)-\frac{r_{kl}}{\ln T}\right|.
\end{eqnarray*}
By the same arguments as those in Lemma 5.3 we obtain $\lim_{T\to \IF} M_{T,2}=0$, and thus the claim follows. \cE{\hfill $\Box$}

\textbf{Proof of Lemma 3.1}:
Using Berman's inequality (\Quan{see for example Piterbarg (1996, Theorem 1.2)}) we have
\begin{eqnarray*}
&&\bigg|P\left\{\max_{t\in \mathfrak{R}(q)\cap\mathbf{S}}X_{k}(t)\leq u_{k}^{(1)}, \max_{t\in\mathfrak{R}(\delta)\cap\mathbf{S}}X_{k}(t)\leq u_{k}^{(2)}, k=1,\cdots,p\right\}\\
&&\ \ \ \ \ \  -P\left\{\max_{t\in \mathfrak{R}(q)\cap \mathbf{S}}\xi_{k}^{T}(t)\leq u_{k}^{(1)}, \max_{t\in\mathfrak{R}(\delta)\cap \mathbf{S}}\xi_{k}^{T}(t)\leq u_{k}^{(2)}, k=1,\cdots,p\right\}\bigg|\\
&&\leq \aa{\mathcal{C}}\sum_{k=1}^{p}\bigg[\sum_{nq\in \mathcal{S}_{i},mq\in \mathcal{S}_{j}\atop i,j=1,2\cdots,n_{T}} \Dkk \int_{0}^{1}\frac{1}{\sqrt{1-r_{kk}^{(h)}(nq,mq)}}\exp\bigg(-\frac{(u_{k}^{(1)})^{2}}{1+r_{kk}^{(h)}(nq,mq)}\bigg)dh\\
&&\ \ \ \ \ \ \ +\sum_{n\delta\in \mathcal{S}_{i},m\delta\in \mathcal{S}_{j}\atop i,j=1,2\cdots,n_{T}} \DkkD \int_{0}^{1}\frac{1}{\sqrt{1-r_{kk}^{(h)}(n\delta,m\delta)}}\exp\bigg(-\frac{(u_{k}^{(2)})^{2}}{1+r_{kk}^{(h)}(n\delta,m\delta)}\bigg)dh\\
&&\ \ \ \ \ \ \ +\sum_{nq\in \mathcal{S}_{i},m\delta\in \mathcal{S}_{j}, nq\neq m\delta\atop i,j=1,2\cdots,n_{T}} \DklQD  \int_{0}^{1}\frac{1}{\sqrt{1-r_{kk}^{(h)}(nq,m\delta)}}\exp\bigg(-\frac{(u_{k}^{(1)})^{2}+(u_{k}^{(2)})^{2}}{2(1+r_{kk}^{(h)}(nq,m\delta))}\bigg)dh\bigg]\\
&&+2\aa{\mathcal{C}}\sum_{1\leq k<l\leq p}\bigg[\sum_{nq\in \mathcal{S}_{i},mq \in \mathcal{S}_{j}\atop i,j=1,2\cdots,n_{T}} \Dkl\int_{0}^{1}\frac{1}{\sqrt{1-r_{kl}^{(h)}(nq,mq)}}\exp\bigg(-\frac{(u_{k}^{(1)})^{2}+(u_{l}^{(1)})^{2}}{2(1+r_{kl}^{(h)}(nq,mq))}\bigg)dh\\
&&\ \ \ \ \ \ \ +\sum_{n\delta\in \mathcal{S}_{i},m\delta\in \mathcal{S}_{j}\atop i,j=1,2\cdots,n_{T}} \DklD \int_{0}^{1}\frac{1}{\sqrt{1-r_{kl}^{(h)}(n\delta,m\delta)}}\exp\bigg(-\frac{(u_{k}^{(2)})^{2}+(u_{l}^{(2)})^{2}}{2(1+r_{kl}^{(h)}(n\delta,m\delta))}\bigg)dh\\
&&\ \ \ \ \ \ \ +\sum_{nq\in \mathcal{S}_{i},m\delta\in \mathcal{S}_{j},nq\neq m\delta\atop i,j=1,2\cdots,n_{T}} \DklQD  \int_{0}^{1}\frac{1}{\sqrt{1-r_{kl}^{(h)}(nq,m\delta)}}\exp\bigg(-\frac{(u_{k}^{(1)})^{2}+(u_{l}^{(2)})^{2}}{2(1+r_{kl}^{(h)}(nq,m\delta))}\bigg)dh\bigg].
\end{eqnarray*}
Now, the claim of Lemma 3.1 \iE{follows from Lemma 5.5-5.6.} \cE{\hfill $\Box$}

{\bf Acknowledgement}: \cH{The authors would like to \iE{deeply}
thank the referees for useful comments and corrections which
improved this paper significantly. E. Hashorva kindly acknowledge
partial support by the Swiss National Science Foundation Grant
200021-1401633/1. Z. Tan's work was \Quan{partially} support by the
National Science Foundation of China Grant 11071182.}


\begin{thebibliography}{100} \small
\bibitem{}
\cE{Adler, R.J.,  {An Introduction to Continuity, Extrema, and Related Topics for General
Gaussian Processes}, Inst. Math. Statist. Lecture Notes Monogr. Ser. 12, Inst. Math.
Statist., Hayward, CA., 1990.}

\bibitem{}
\cE{Berman, M.S., Limit theorems for the maximum term in stationary
sequences.  Ann. Math. Statist., 1964, 35,  502-516.}


\bibitem{} Berman, S.M., Sojourns and extremes of Gaussian processes, Ann. Probab., 1974, 2, 999-1026.

\bibitem{} \cT{Berman,  M.S., Sojourns and Extremes of Stochastic Processes. Wadsworth \& Brooks/ Cole, Boston, 1992.}


\bibitem{17} {D\c{e}bicki, K., Ruin probability for Gaussian integrated processes. Stoch. Proc. Appl., 2002,   98, 151-174.}


\bibitem{}
D\c{e}bicki, K.,  Michna, Z., Rolski, T., Simulation of the asymptotic constant in some fluid models. Stochastic Models, 2003, 19, 407-423.

\bibitem{}
\cE{D\c{e}bicki, K., Kisowski, P., A note on upper estimates for Pickands constants. Stat. Prob. Letters, 2009, 78, 2046-2051.}

\bibitem{DTA} \AA{D\c{e}bicki, K., Tabi\'{s}, K.,  Extremes of time-average stationary Gaussian processes. Stoch. Proc. Appl., 2011, 121,
    2049-2063.}


\bibitem{} \ccE{Hashorva, E., Ji, L., Piterbarg, V.I., On the supremum of gamma-reflected processes with fractional Brownian motion as input.
 Stoch. Proc. Appl., 2013, 123,11, 4111-4127.}

\bibitem{} \ccE{Hashorva, E., Peng, Z., Weng, Z.,
On Piterbarg theorem for the maxima of stationary Gaussian sequences. Lithuanian Math. J., 2013, 53,3, 280-292.}



\bibitem{} Hashorva, E., Weng, Z., Limit laws for extremes of dependent
stationary Gaussian arrays. Stat. Probab. Letters, 2013, 83, 320-330.


\bibitem{} H\"{u}sler, J., Dependence between extreme values of discrete and continuous time locally stationary Gaussian processes.
Extremes, 2004, 7, 179-190.

\bibitem{} H\"{u}sler, J.,  Sch\"{u}pbach, M., Limit results for maxima in nonstationary multivariate Gaussian sequences. Stoch. Proc. Appl.,
1988, 27, 91-99.

\bibitem{} H\"{u}sler, J., Piterbarg, V.I., Limit theorem for maximum of the storage process with fractional Brownian motion as input.
Stoch. Proc. Appl., 2004, 114, 231-250.

\bibitem{}James, B., James, K. and Qi, Y.  Limit distribution of the sum and maximum from multivariate Gaussian sequences.
J. Multi. Analysis. 2007, 98, 517-532


\bibitem{} Leadbetter, M.R., Lindgren, G. and Rootz\'{e}n, H., Extremes and Related Properties of Random Sequences and
Processes. Series in Statistics, Springer, New York, 1983.

\bibitem{} McCormick, W. P. and Qi, Y. Asymptotic distribution for the
sum and the maximum of Gaussian processes. J. Appl. Probab. 2000, 37, 958-971.

\bibitem{}
Mishura, Y. and Valkeila, E.  An extension of the L\'{e}vy characterization to fractional Brownian motion. Ann. Probab., 2011, 39, 439–-470.

\bibitem{} Mittal, Y. and Ylvisaker, D. Limit distribution for the
maximum of stationary Gaussian processes. Stoch. Proc. Appl., 1975, 3, 1-18.


\bibitem{} Pickands, J., III, Asymptotic properties of the maximum in a stationary Gaussian process. Trans. Am. Math. Soc., 1969, 145, 75-86.

\bibitem{} \cE{
Piterbarg, V.I.,  On the paper by J. Pickands "Upcrosssing probabilities for stationary Gaussian processes". Vestnik Moscow. Univ. Ser. I Mat. Mekh. 27, 25-30. English transl. in Moscow Univ. Math. Bull., 1972, 27.}

\bibitem{} Piterbarg, V.I., Asymptotic Methods in the Theory of Gaussian Processes and Fields, AMS, Providence, 1996.

\bibitem{} Piterbarg, V.I., Discrete and continuous time extremes of Gaussian processes, Extremes, 2004, 7, 161-177.


\bibitem{} \cE{Tan, Z., Hashorva, E., On Piterbarg max-discretisation theorem for standardised maximum of stationary
Gaussian processes. Meth. Computing Appl. Probab., 2014, 16,1, 169-185.}


\bibitem{} \cE{Tan, Z., Hashorva, E., Peng, Z., Asymptotics of maxima of strongly dependent Gaussian processes. J. Appl. Probab. 2012, 49, 1106-1118.}

\bibitem{} Tan, Z., Tang, L., The dependence of extremes values of discrete and continuous time strongly dependent Gaussian
processes. Stochastics, 2012,  DOI: 10.1080/17442508.2012.756489.


\bibitem{TW12} Tan, Z., Wang, Y.,  Some asymptotic results on extremes of incomplete samples.
Extremes, 2012, 15, 319-332.


\bibitem{TW11} \cE{Tan, Z., Wang, Y.,  Extreme values of discrete and continuous time strongly dependent Gaussian processes.
Comm. Stat. Theory Meth., 2013, 42,  2451-2463,}


\bibitem{} \iE{Turkman, K.F., Discrete and Continuous time series extremes of stationary processes.
Handbook of statistics  Vol 30.  Time Series Methods and Aplications. Eds.
T.S. Rao, S.S. Rao and C.R. Rao. Elsevier, pp. 565-580, 2012.}

\bibitem{}
Wu, D. Generalized Pickands constants. Journal of Mathematical Physics, 2007, 48, 053513-053513.

\end{thebibliography}
\end{document}